\documentclass[11pt]{amsart}

\usepackage[english]{babel}
\usepackage{graphicx,amsmath,amssymb,color, enumerate}
\numberwithin{equation}{section}
\usepackage{enumerate}
 \usepackage{xcolor}
 \definecolor{brickred}{rgb}{0.8, 0.25, 0.33}
\definecolor{blue(ryb)}{rgb}{0.01, 0.28, 1.0}
\definecolor{brandeisblue}{rgb}{0.0, 0.44, 1.0}
\definecolor{ceruleanblue}{rgb}{0.16, 0.32, 0.75}
\definecolor{cobalt}{rgb}{0.0, 0.28, 0.67}
\definecolor{coolblack}{rgb}{0.0, 0.18, 0.39}
\definecolor{darkblue}{rgb}{0.0, 0.0, 0.55}
\usepackage[hidelinks]{hyperref}
\hypersetup{
    colorlinks,
    citecolor=darkblue,
    filecolor=black,
    linkcolor=darkblue,
    urlcolor=black
}

\newtheorem{theorem}{Theorem}[section]
\newtheorem{lemma}[theorem]{Lemma}
\newtheorem{proposition}[theorem]{Proposition}

\theoremstyle{remark}

\newtheorem{remark}[theorem]{Remark}

\newcommand{\sq}{\mathsf{q}}
\newcommand{\sw}{\mathrm{sw}}
\newcommand{\dw}{\mathrm{dw}}
\newcommand{\wkb}{\mathrm{WKB}}
\newcommand{\harm}{\mathrm{harm}}

\newcommand{\ee}{\mathrm{e}}
\newcommand{\dd}{\mathrm{d}}

\newcommand{\cO}{\mathcal{O}}
\newcommand{\cS}{\mathcal{S}}

\newcommand{\N}{\mathbb{N}}
\newcommand{\Z}{\mathbb{Z}}
\newcommand{\R}{\mathbb{R}}

\newcommand{\Ab}{\mathbf{A}}
\newcommand{\Fb}{\mathbf{F}}
\newcommand{\cB}{\mathcal B}
\newcommand{\cD}{\mathcal D}
\newcommand{\cE}{\mathcal E}
\newcommand{\cK}{\mathcal K}

\newcommand{\cJ}{\mathcal J}
\newcommand{\cL}{\mathcal L}
\newcommand{\cR}{\mathcal R}

\newcommand{\sI}{\mathsf I}

\newcommand{\sH}{\mathsf{H}}
\newcommand{\sT}{\mathsf{T}}
\newcommand{\sW}{\mathsf{W}}
\newcommand{\sU}{\mathsf{U}}

\newcommand{\sR}{\mathsf{R}}
\newcommand{\sL}{\mathsf{L}}
\newcommand{\sK}{\mathsf{K}}

\newcommand{\ii}{\,\mathrm{i}}
\newcommand{\xv}{\mathsf{x}}
\newcommand{\yv}{\mathsf{y}}

\DeclareMathOperator{\curl}{curl}
\newcommand{\cT}{\mathcal{T}}
\newcommand{\vb}{\mbox{\emph{\textbf{v}}}}
\newcommand{\bp}{\mathbf{p}}
\newcommand{\bq}{\mathbf{q}}

\begin{document}

\title[Quantum Tunneling and the Aharonov-Bohm effect]{Quantum Tunneling and the Aharonov-Bohm effect}

\author{Bernard Helffer}
\address[B. Helffer]{Laboratoire de Math\'ematiques Jean Leray, CNRS, Nantes Universit\'e,
44000 Nantes, France.}
\email{Bernard.Helffer@univ-nantes.fr}

\author{Ayman Kachmar}
\address[A. Kachmar]{School of Science and Engineering, The Chinese University of Hong Kong Shenzhen, Guangdong, 518172, P.R. China.}
\email{akachmar@cuhk.edu.cn}

\date{\today}

\begin{abstract}
We investigate a Hamiltonian with radial potential wells and an Aharonov-Bohm vector potential with two poles. Assuming that the potential wells are symmetric, we derive the semi-classical asymptotics of the splitting between the ground and second state energies. The  flux effects due to the Aharonov-Bohm vector potential are of lower order compared to the contributions coming from the potential wells.
\end{abstract}
\maketitle

\section{Introduction}

\subsection{Motivation}

In the presence of symmetric potential wells, the eigenvalue splitting between the ground and second state energies of the semi-classical Schr\"odinger operator is exponentially small, and the exponential decay involves the (Agmon) distance between the wells \cite{HSj, Si}. This phenomenon, known as quantum tunneling in the literature, can be seen as a result of the structure of the eigenfunctions, which are in the form of a superposition of functions localized near the potential wells, and thanks to symmetry, they yield equal probabilities for the particle to be in either well. 

Under a magnetic field, the eigenfunctions are no more real-valued and have non-trivial phases, making their approximation and  decay more subtle to capture.  Under a uniform  magnetic field and symmetric radial potential wells, progress in this direction was the subject of the recent works \cite{FSW, HK, Mo}, and for instance one can observe a change of multiplicity and eigenvalue crossings as a  flux and symmetry effect  \cite{HKS}. 

If no external potential is present, quantum tunneling under uniform magnetic field  can be induced by the geometry  \cite{BHR}, and following a similar proof, by a discontinuity and sign change of the magnetic field \cite{FHK, Abal}. Quite recently, a first example of a  quantum tunneling  as a result of a non-homogeneous magnetic field with symmetric radial magnetic wells was established in \cite{FMR}.

In their celebrated paper \cite{AB},  Aharonov and Bohm introduced a vector potential generating a zero-magnetic field except in a region with `zero' radius,  while the total flux remains fixed. Despite that there is no magnetic field, the charged particle `feels' the magnetic flux, for instance through the `energy levels', and this phenomenon is known as the Aharonov-Bohm effect. It has motivated several mathematical works, for instance \cite{AT, CF, H1, HHOO}.   

In this paper, our aim   is to understand the Aharonov-Bohm effect on quantum tunneling.  Interest in this question in physics traces back  to a  paper published in 2014 on the  tunneling rate of certain structures \cite{NST}, where a periodic dependence on the  magnetic flux was reported. In the setting of symmetric potential wells, we study how does the (magnetic) vector potential influence the eigenvalue splitting and the structure of the eigenfunctions.  Our main results in Theorems~\ref{thm:main} and \ref{thm:main*} quantify the dependence of the eigenvalue splitting (which corresponds to the magnitude of tunneling)  on the magnetic flux.

\subsection{The Aharonov-Bohm potential}

We introduce the vector potential
\begin{equation}\label{eq:def-AB}
\R^2\ni\xv=(x_1,x_2)\mapsto\Fb(\xv)=\begin{pmatrix}-\displaystyle\frac{x_2}{|\xv|^2}  , \frac{ x_1}{|\xv|^2}\end{pmatrix},
\end{equation}
and we notice that it generates zero magnetic field everywhere except at the origin of $\R^2$. In fact, as a distribution, the magnetic field is a multiple of the Dirac delta function supported at the origin:
\[\curl\Fb=2\pi \, \delta_0\quad \mbox{in }\cD'(\R^2). \] 
The circulation of $\Fb$ along a loop $C$  around the origin is always non-zero and independent of the loop,
\[\frac1{2\pi}\int_{C}\Fb\cdot\,\dd r=1,\]
and we interpret it as the magnetic flux induced by $\Fb$. Observing that $\Fb$ induces unit flux, it is natural to consider $\alpha>0$ and the vector potential
\begin{equation}\label{eq:def-F}
\Fb_\alpha (\xv)= \alpha \Fb(\xv),
\end{equation}  
and we call $\Fb_\alpha$ the Aharonov-Bohm potential with flux $\alpha$ and pole $0$.

There is a Hardy type inequality \cite{LW} that shows the flux effect through  the following $1$-periodic function
\begin{equation}\label{eq:def-e}
e(\cdot)=\inf_{m\in\Z}|\cdot-m|,
\end{equation}
and which reads for all $h,\alpha>0$ and $u\in C_c^\infty(\R^2\setminus\{0\})$ as
\begin{equation}\label{eq:Hardy}
\int_{\R^2}|(-\ii h\nabla-\Fb_\alpha)u|^2\dd\xv \geq h^2e(\alpha/h)\int_{\R^2}\frac{|u|^2}{|\xv |^2}\dd\xv .
\end{equation}
It can be easily observed by expressing $u$ in polar coordinates $(r,\theta)$ as
\[ u=\sum_{m\in\Z}u_m(r)\ee^{\ii m\theta},\]
and by noticing  that its `kinetic energy' under the potential $\Fb_\alpha$ is
\[
\int_{\R^2}|(-\ii h\nabla-\Fb_\alpha)u|^2\dd\xv =2\pi\sum_{m\in\Z}\int_{\R_+}\Bigl(h^2|u_m'(r)|^2+ \frac{(mh-\alpha)^2}{r^2} |u_m(r)|^2\Bigr)r\,\dd r.
\]

\subsection{Single potential well}

Suppose that $\vb \in C_c^\infty(\mathbb R^2)$ is a  real-valued 
radially symmetric function with a unique non-degenerate minimum at $0$.   More precisely, we assume that $\vb(\xv )=v(|\xv |)$, where $v$ belongs to $C^\infty(\overline{\mathbb R}_+)$ and, for some $\sigma >0$, 
\begin{equation}\label{eq:cond-v}
{\rm supp}\,v\subset [0,\sigma],\quad k:=\min v<0,\quad v^{-1}(k)=0,\quad v''(0)>0.
\end{equation}
Notice that, since $\vb \in C^\infty_c(\R^2)$ is radial, we have  $\vb (r,0)=\vb (-r,0)$ for any $r\in\mathbb R$,  which yields that $\frac{\dd^n}{\,\dd r^n} v(r)|_{r=0}=0$ for all odd integers $n$. Hence we can extend $v$ to an even $C^\infty$- function on $\mathbb R$.

With $\alpha,h>0$ and  $\Fb_\alpha$ introduced in \eqref{eq:def-F}, consider the self-adjoint operator\footnote{We use the convention $(-\ii h\nabla-\Fb_\alpha)^2=(-\ii h\nabla-\Fb_\alpha)\cdot(-\ii h\nabla-\Fb_\alpha)$, and since $\mathrm{div}\Fb_\alpha=0$, we get $-h^2\Delta+2\ii h\Fb_\alpha\cdot\nabla+|\Fb_\alpha|^2$ on $\R^2\setminus\{0\}$.} in $L^2(\R^2)$,
\begin{equation}\label{eq:def-H0}
\sH_{\alpha,0}:=(-\ii h\nabla-\Fb_\alpha)^2+\vb. 
\end{equation}
which is defined as the  Friedrichs extension, starting from the quadratic form
\begin{equation}\label{eq:def-q-alpha}
C_c^\infty(\R^2\setminus\{0\})\ni u\mapsto \sq_\alpha(u):=\int_{\R^2}\bigl(|(-\ii h\nabla-\Fb_\alpha)u|^2+\vb |u|^2\bigr) \dd\xv .
\end{equation}
In the same way, we get a self-adjoint realization of $P_\alpha:=(-\ii h\nabla-\Fb_\alpha)^2$ in $L^2(\R^2)$. When $v$ satisfies \eqref{eq:cond-v}, it  is $P_\alpha$-compact, and  $\sH_{\alpha,0}$ has the same domain and essential spectrum as $P_\alpha$. Moreover, for $h$ sufficiently small, the ground state energy  of $\sH_{\alpha,0}$ is negative and therefore belongs to the discrete spectrum.

Throughout this paper, we fix $\alpha>0$ and $e_0\in[0,\frac12]$, and suppose that $h$ varies in the following set
\begin{equation}\label{eq:def-J}
\mathcal J_\alpha(e_0)=\{h\in\R_+\colon e(\alpha/h)=e_0\},
\end{equation}
where $e(\cdot)$ is introduced in \eqref{eq:def-e}. 

Let us point out a few remarks  regarding the operator $\sH_{\alpha,0}$:
\begin{enumerate}[--]
\item We exclude $\alpha=0$ because this corresponds to the pure electric operator $\sH_{0,0}:=-h^2\Delta+\vb $.
\item It suffices to consider the case of $\alpha>0$ because the operators $\sH_{\alpha,0}$ and $\sH_{-\alpha,0}$ are unitarily equivalent (by applying the unitary transformation of complex conjugation).
\item  If $h\in\cJ_\alpha(e_0)$ and $e_0=0$, then $\alpha/h$ is an integer  and the operators $\sH_{\alpha,0}$ and $\sH_{0,0}$ are unitarily equivalent, thanks to the following identity 
\[ \ee^{-\ii\frac{\alpha}{h}\theta}(-\ii h\nabla-\Fb_\alpha)\ee^{\ii\frac{\alpha}{h}\theta}=-\ii h\nabla.\]
\item The operators $\sH_{\alpha,0}$ and $\sH_{\hat\alpha,0}$, with $\hat\alpha=\alpha+1/h$, are unitarily equivalent.
\end{enumerate}

\medskip The structure of the ground states of $\sH_{\alpha,0}$ involves the following (Agmon) distance to the origin
\begin{equation}\label{eq:dist}
d(r)=\int_0^r\sqrt{v(\rho)-v(0)}\,\dd \rho \qquad(r\geq 0),
\end{equation}
and the following function
\begin{equation}\label{eq:def-phase}
p_{e_0}(r):=\frac{d''(r)+(1+2e_0)\frac{d'(r)}{r}-\sqrt{2v''(0)}(1+e_0)}{2d'(r)}\qquad(r\geq 0).
\end{equation}
Thanks to \eqref{eq:cond-v}, we have $p_{e_0}(r)=\cO(1)$ as $r\to0$, hence the  function defined as
\begin{equation}\label{eq:def-hat-a}
\begin{gathered}\widehat a_0(r)=\widehat a _0(r,e_0):=A_0\exp\Bigl(-\int_0^rp_{e_0}(\rho)\dd\rho\Bigr),\\
A_0=A_0(v,e_0):= 2^{(1-e_0)/4}\,v''(0)^{(1+e_0)/4}\, \Gamma(1+e_0)^{-1/2}\,,
\end{gathered}
 \end{equation}
is continuous on $\overline\R_+$.

We are ready now to describe the flux effects on the approximation of the ground states of $\sH_{\alpha,0}$, which will distinguish between the case where $\alpha/h$ is not a half-integer ($0\leq e_0<\frac12$), and the case of $\alpha/h$ being a half-integer  ($e_0=\frac12$). 

For $t \in\R$, we denote by $m_*(t)$ the smallest integer satisfying \[e(m_*(t))=|t-m_*(t)|\,.\]
 Note that,  $e_*(t)=|t-m|$ has a unique integer solution when $t\not\in\Z+\frac12$, while there are two integer solutions $m_*(t),m_*(t)+1$ in the case where $t\in\Z+\frac12$.
\begin{theorem}\label{thm:sw}
Suppose that $\alpha>0$ and $0\leq e_0\leq 1/2$, and let the set $\cJ_\alpha(e_0)$ be as in \eqref{eq:def-J}.
\begin{enumerate}[\rm i)]
\item If $0\leq e_0<\frac12$, then there is $h_0>0$ such that, if $h\in \cJ_\alpha(e_0)\cap (0,h_0]$, the ground state energy of $\sH_{\alpha,0}$ is a simple eigenvalue, and it has a normalized ground state $\phi$ that satisfies locally uniformly on $\R^2\setminus\{0\}$,
\[ \Bigl|\ee^{d(r)/h}\phi(r,\theta)-\pi^{-\frac12}h^{-\frac{1+e_0}{2}}r^{e_0}\widehat a_0(r,e_0)\ee^{\ii m_*(\alpha/h)\theta}\Bigr|=\cO\bigl( h^{\frac{1-e_0}{2}}\bigr), \]
where $(r,\theta)$ are the polar coordinates in $\R^2$, and $d(r)$ is as in \eqref{eq:dist}.
\item   If $e_0=\frac12$, then there is $h_1>0$ such that, if $h\in\cJ_\alpha(e_0)\cap(0,h_1]$, the ground state energy of $\sH_{\alpha,0}$ has multiplicity $2$. Moreover, there are two normalized ground states
$\phi_1$ and $\phi_2$ such that
\[\phi_2=\ee^{\ii \theta}\phi_1\]
and $\phi_1$ satisfies  locally uniformly on $\R^2\setminus\{0\}$,
\[
 \Bigl|\ee^{d(r)/h}\phi_1(r,\theta)-\pi^{-\frac12}h^{-\frac{3}{4}}r^{e_0}\widehat a_0(r,\mbox{$\frac12$})\ee^{\ii m_*(\alpha/h)\theta}\Bigr|=\cO\bigl( h^{\frac{1}{4}}\bigr).\\
\]
\end{enumerate}
\end{theorem}
If we consider the vector potential $\Fb_\alpha(\xv -\xv_0)$, with $\xv_0\in\R^2$,  the ground states are still localized near the minimum of $\vb$. This minimum being $0$, the  effect of the potential $\Fb_\alpha(\xv -\xv_0)$ will disappear  when $\xv_0\not=0$. Consequently, to observe the flux effect (to leading order) in the semi-classical limit, the pole of the vector potential has to be the same as the potential well.

\subsection{Double potential wells and tunneling}~\medskip

We now consider as in \cite{FSW}  the case of a potential with two wells,
\begin{equation}\label{eq:def-V}
V(\xv )=v(|\xv -\xv_\ell|)+ v(|\xv -\xv_r|),
\end{equation}
where 
\[\xv_\ell=(-L/2,0),\quad \xv_r=(L/2,0),\quad L>2\sigma ,\]
and $v$ is the function introduced in \eqref{eq:cond-v}.

Moreover, we consider the vector potential
\begin{equation}\label{eq:def-A}
\Ab_\alpha (\xv )=\Fb_\alpha(\xv -\xv_\ell)+  \Fb_\alpha(\xv -\xv_r),
\end{equation}
 with two poles centered at the wells $\xv_\ell$ and $\xv_r$, where $\alpha>0$ and $\Fb_\alpha$ is the Aharonov-Bohm potential introduced in \eqref{eq:def-F}.

Consider the magnetic Schr\"odinger operator in $L^2(\R^2)$,
\begin{equation}\label{eq:def-op}
\sH_\alpha=(-\ii h\nabla-\Ab_\alpha)^2+V, 
\end{equation}
and the associated quadratic  form 
\[u\mapsto \int_{\R^2}\bigl(|(-\ii h\nabla-\Ab_\alpha)u|^2+V(\xv )|u|^2\bigr)\dd\xv .\]
Initially defined on $C_c^\infty(\R^2\setminus\{\xv_\ell,\xv_r\})$, this quadratic form  is semi-bounded, and with 
form domain 
 the magnetic Sobolev space,
\[H^1_{\Ab_\alpha,h}(\R^2):=\{u\in L^2(\R^2)\colon (-\ii h\nabla -\Ab_{\alpha})u\in L^2(\R^2)\},\]
it is moreover closed, hence  by the Friedrichs theorem, it realizes $\sH_\alpha$ as a self-adjoint operator in $L^2(\R^2)$.  

If $\alpha\in h\Z$,  then  $\sH$ is unitarily equivalent to $\sH_0:=-h^2\Delta+V$ \cite{HHOO}. Whereas, if $\alpha\not\in h\Z$,  it results from the Hardy inequality in \eqref{eq:Hardy} that the  form domain is (see also \cite[Sec. 2 and Rem. 3.8]{KP})
\[H^1(\R^2)\cap L^2(\R^2;|\xv -\xv_\ell|^{-2}\dd\xv )\cap L^2(\R^2;|\xv -\xv_r|^{-2}\dd\xv ).\] 
We denote by $(\lambda_j(h,\alpha))_{j\geq 1}$ the  sequence of min-max eigenvalues of $\sH$, and we assume that $h$ varies in the set $\cJ_\alpha(e_0)$ introduced in \eqref{eq:def-J}.
We  study first the case where
\[ \alpha/h\not\in\Z+\frac12\,,
\]
and obtain the following flux dependence of the eigenvalue splitting.
\begin{theorem}\label{thm:main}
Suppose that $ L>2\sigma$, $\alpha>0$ and $e_0\in[0,\frac12)$.  Let
\[S(v,L)=2\int_0^{L/2}\sqrt{v(r)-v(0)}\, \,\dd r. \]
Then, there exist an explicit  constant $C(L,v,e_0)>0$ and a function $f_{e_0}:\cJ_\alpha(e_0)\to\R$ such that, $f_{e_0}(h)\to0$ as $h\to0$, and, for $h\in\cJ_\alpha(e_0)$, we have
\[\lambda_2(h,\alpha)-\lambda_1(h,\alpha) =
C(L,v,e_0) \bigl(1+ f_{e_0}(h)\bigr) h^{\frac12-e_0} \ee^{-S(v,L)/h}.  \]
Moreover, for any positive $\delta_0<\frac12$, $f_{e_0}(h)\to 0$ uniformly with respect to $e_0\in [0,\delta_0]$.
\end{theorem}
\begin{remark}~
\begin{enumerate}[i)]
\item  The eigenvalue splitting involves the Agmon distance $S(v,L)$ between the two wells associated with the double well  potential $V-v(0)$,  while the flux appears in the prefactor, which amounts to an Aharonov-Bohm effect on the tunneling between two potential wells.
\item The definition of the constant $C(L,v,e_0)$  involves the function $p_{e_0}$ introduced in \eqref{eq:def-phase}:
\begin{equation}\label{eq:CLv}
C(L,v,e_0)=\frac{2^{\frac{4-5e_0}2}|v''(0)|^\frac{1+e_0}{2}L^{2e_0+\frac12}}{\pi^{\frac12}|v(0)|^{\frac14}\Gamma(1+e_0)}\exp\Bigl(-2\int_0^{L/2}p_{e_0}(\rho)\dd\rho \Bigr)\,.
\end{equation} 
\item  If   $e_0=0$, which corresponds to $\alpha/h\in\Z$, then the flux effects disappear. The tunneling formula  
\[\lambda_2(h,0)-\lambda_1(h,0)\underset{h\to0}{=} 
C(L,v,0) (1+o(1))\, h^{\frac12} \ee^{-S(v,L)/h}\]
recovers a particular case of the general double well analysis of  Helffer-Sj\"ostrand \cite{HSj} and Simon \cite{Si}  for the Schr\"odinger operator $-h^2\Delta+V$. This is in fact a consequence of the unitary equivalence of $\sH_\alpha$ and $\sH_0$. 
\item  Along the proof of Theorem~\ref{thm:main}, we get $\lambda_3(h,\alpha)\geq h\sqrt{v''(0)/2}+o(h)$, and if $0<e_0<\frac12$,
\begin{equation}\label{eq:lambda4-3}
\lambda_4(h,\alpha)-\lambda_3(h,\alpha) =C(L,v,1-e_0)h^{e_0-\frac12}\bigl(1+o(1)\bigr)  \ee^{-S(v,L)/h},
\end{equation}
which interestingly shows a change in the order of the prefactor when compared with the asymptotics of 
$\lambda_2(h,\alpha)-\lambda_1(h,\alpha)$. 
\end{enumerate}
\end{remark}
\medskip
 A similar result  holds in the case where $e_0=\frac12$, which follows by computing the eigenvalues of a $4\times4$ interaction matrix.   To leading order, the quantity $C(L,v,\mbox{$\frac12$})  \ee^{-S(v,L)/h}$ will 
 give the scale for  the difference between two consecutive eigenvalues. 

\begin{theorem}\label{thm:main*}
Suppose that $ L>2\sigma$, $\alpha>0$ and $e_0=\frac12$. Then, 
there exists $h_0>0$ such that, for  $h\in(0,h_0]\cap\cJ_{e_0}(\alpha)$, it holds for the four lowest eigenvalues
\[\begin{aligned}
\lambda_2(h,\alpha)-\lambda_1(h,\alpha) &=C(L,v,\mbox{$\frac12$})  \ee^{-S(v,L)/h}+o\bigl(\ee^{-S(v,L)/h} \bigr),\\
\lambda_3(h,\alpha)-\lambda_2(h,\alpha) &=o\bigl(\ee^{-S(v,L)/h} \bigr),\\
\lambda_4(h,\alpha)-\lambda_3(h,\alpha) &=C(L,v,\mbox{$\frac12$})  \ee^{-S(v,L)/h}+o\bigl(\ee^{-S(v,L)/h} \bigr),\\
\lambda_5(h,\alpha)-\lambda_4(h,\alpha)&\geq 3h\sqrt{v''(0)/2}+o(h) .
\end{aligned}
\]
\end{theorem}
  Theorem~\ref{thm:main*} leaves open   the possibility of  crossings  or equality between the second and the third eigenvalues for the double well operator, but we give  indications in Subsection~\ref{subsec:mult} where multiplicity could occur as a consequence of the symmetries of the problem.

\subsection{Organization} The rest of the paper consists of six sections devoted to the proof of the theorems announced in the introduction. In Section~\ref{sec:1D-op}, we study the ground state energy and the accurate approximation of the ground states for a one dimensional operator, which we use later in Section~\ref{sec:sw-op}  to study the single well operator in \eqref{eq:def-H0} and to prove Theorem~\ref{thm:sw}. To prove Theorem~\ref{thm:main},  we construct in Section~\ref{sec:qm} quasi-modes for the double well operator in \eqref{eq:def-op},  enabling us to reduce to a $2\times2$ matrix, then we calculate the eigenvalues of this matrix in Section~\ref{sec:interaction}. Finally, we discuss the case of the   half-integer flux in Section~\ref{sec:interaction*} and prove Theorem~\ref{thm:main*}.

\section{Study of a one dimensional operator}\label{sec:1D-op}

 Throughout this section, we fix  $e_0\in[0,\frac12]$, and for $h>0$, we consider the operator
\begin{equation}\label{eq:def-T0}
\sT=\sT_{e_0}:=-h^2\frac{\dd^2}{\,\dd r^2}-\frac{h^2}{r}\frac{\dd}{\,\dd r}+v(r) +\frac{h^2e_0^2}{r^2},
\end{equation}
which will play a crucial role in  proving Theorem~\ref{thm:sw}. We are interested 
in the case where $v\in C^\infty(\overline{\R}_+)$ and satisfies \eqref{eq:cond-v}, but we will encounter also the harmonic potential where  $v(r)=\beta r^2$  and $\beta>0$.
In both cases, we see that $v$ is smooth and bounded from below.

\subsection{Domain and structure of ground states}
The operator $\sT$ is  the  self-adjoint operator in $L^2(\R_+,r\,\dd r)$ corresponding to the  quadratic form
\begin{equation}\label{eq:def-form-T0}
\sq(u)=\int_{\R_+}\Bigl(h^2|u'(r)|^2+v(r)|u(r)|^2+\frac{h^2e_0^2}{r^2}|u(r)|^2\Bigr)r\,\dd r,
\end{equation}
initially defined on $C_c^\infty(\R_+)$. With $C>-\min v$, the closure of $C_c^\infty(\R_+)$ under the norm $(C\|u\|^2_{L^2(\R_+,r\,\dd r)}+\sq(u))^{1/2}$ is
\begin{equation}\label{eq:def-form-dom-T0}
\mathrm{Dom}(\sq)=\{u\in L^2(\R_+,r\,\dd r)\colon u', u/r\in L^2(\R_+,r\,\dd r) \},
\end{equation}
and the domain of $\sT$ is
\begin{equation}\label{eq:dom-T0}
\mathrm{Dom}(\sT)=\{u\in \mathrm{Dom}(\sq)\colon \sT u\in L^2(\R_+,r\,\dd r)\}. 
\end{equation} 
In particular, if $v$ satisfies \eqref{eq:cond-v}, $u$ belong to $\mathrm{Dom}(\sT)$ whenever
\[u,u',u/r,-u''-u'/r+e_0^2u/r^2\in L^2(\R_+,r\,\dd r). \]

We denote by $\{\lambda_n(\sT)\}_{n\geq 1}$ the  min-max spectral sequence corresponding to $\sT$, and note that if $v$ satisfies \eqref{eq:cond-v}, then the essential spectrum of $\sT$ is $[0,+\infty)$.

Our aim is to prove the following theorem.
\begin{theorem}\label{thm:T0}
Suppose that $v$ satisfies \eqref{eq:cond-v}. Then, there is $h_0>0$ such that, for every $h\in (0,h_0]$, $\lambda_1(\sT)$ is a simple eigenvalue, with positive normalized ground state $\psi$ that satisfies locally uniformly on $\R_+$,
\[
\begin{gathered}
 \Bigl|\ee^{d(r)/h}\psi(r)-h^{-\frac{1+e_0}{2}}r^{e_0}\widehat a_0(r,e_0)\Bigr|=\cO\bigl( h^{\frac{1-e_0}{2}}\bigr),\\
  \Bigl|\ee^{d(r)/h}\psi'(r)+h^{-\frac{3+e_0}{2}}d'(r)r^{e_0}\widehat a_0(r,e_0)\Bigr|=\cO\bigl( h^{-\frac{2+e_0}{2}}\bigr),
  \end{gathered} \]
where $d(r)$ and $\widehat a_0(r,e_0)$ are  as in \eqref{eq:dist} and \eqref{eq:def-hat-a}, respectively.
\end{theorem}
\begin{remark}\label{rem:T0-unif}
In the proof of Theorem~\ref{thm:T0}, we can follow the dependence on $e_0$ in the remainder terms.  We will see that that the $e_0$-dependent contributions in the remainder are  controlled by sum of  terms of the form  $a_n e_0^n c_n ^{e_0}$. 
 Therefore,
 the estimate in Theorem~\ref{thm:T0}, as well as all the estimates written in this section, hold uniformly with respect to $e_0\in[0,\frac12]$.
\end{remark}

\begin{remark}\label{rem:T00}
In the case $e_0=0$, $\sT$ corresponds to the radial part of\break$\sH=-h^2\Delta+\vb $, the Laplacian with a radial potential, and the conclusion of Theorem~\ref{thm:T0} is known.
\end{remark}

\begin{remark}\label{rem:T0-flat}
Thanks to the unitary transformation 
\[U\colon L^2(\R_+,\,\dd r)\ni f\mapsto \sqrt{r}f \in L^2(\R_+,r\,\dd r),\] the operator $\sT$ is unitarily equivalent to\footnote{
The domain of $\widetilde T$ is $\mathrm{Dom}(\widetilde\sT)=\{f\in H^1_0(\R_+)\colon \widetilde\sT f\in L^2(\R_+)\}$.}
\[\widetilde \sT=-h^2\frac{\dd^2}{\,\dd r^2}+\frac{e_0^2-1/4}{r^2}+ v. \]
In particular, when $e_0=1/2$, the operator $\widetilde\sT$ becomes the one dimensional Schr\"odinger operator with a smooth potential  and with Dirichlet boundary condition\footnote{ If $u$ is in the form domain of $\widetilde\sT$ and $u(0)\not=0$, then by continuity of $u$, $(e_0/r)u\not\in L^2(\R_+,r\,\dd r)$.}  on $\R_+$, but with potential minimum at $0$.
\end{remark}

\subsection{Harmonic potential}

In the case where $v(r)=\beta r^2$, with $\beta>0$, we can reduce to the case $h=1$ by scaling, and deal with the operator
\[\sT^{\harm}=-\frac{\dd^2}{\,\dd r^2}-\frac1r\frac{\dd}{\,\dd r}+\beta r^2+\frac{e_0^2}{r^2}.\]
The spectrum of $\sT^{\harm}$ is purely discrete,  and  consists of the simple eigenvalues  \cite[Sec.~2]{H},
\begin{equation}\label{eq:2.2}
E_n(e_0)=2\sqrt{\beta}(1+e_0+2n),\quad (n=0,1,2,\cdots).
\end{equation}
In particular, we find the  ground state energy 
\begin{equation}\label{defmu1}
E_0(e_0)=2\sqrt{\beta}\bigl(1+e_0\bigr)
\end{equation}
with the positive  and normalized ground state
\begin{equation}\label{eq:def-u1}
u^\harm(r)= 2^{\frac 12} \beta^{\frac{1+e_0}{4}}\, \Gamma(e_0+1)^{-\frac 12}\,r^{e_0}\ee^{-\sqrt{\beta}r^2/2}.
\end{equation}
In the case $e_0=0$, we recover in \eqref{eq:2.2} the eigenvalues with  radial eigenfunctions  for the
two dimensional harmonic oscillator,  $-\Delta +\beta|\xv |^2$.

\subsection{Harmonic approximation}

In general, if we assume that $v$ satisfies \eqref{eq:cond-v}, we can use the harmonic approximation, and obtain that the eigenvalues are  effectively close to those corresponding to the harmonic potential $\vb (\xv )=\beta|\xv |^2$, with $\beta=v''(0)/2$. 
\begin{proposition}\label{prop:harm-app}
Suppose that $v$ satisfies \eqref{eq:cond-v}, then the ground state energy of $\sT$ satisfies
\[\lambda_1(\sT)=v(0)
+h\sqrt{2v''(0)}\bigl(1+e_0\bigr)+o(h) \mbox{ as }h\to0. \]
Moreover, 
\[\lambda_2(\sT)-\lambda_1(\sT)\geq 
h\sqrt{2v''(0)} +o(h).\]
\end{proposition}
Since $v(0)<0$ and since the essential spectrum of $\sT$ is equal to $[0,+\infty)$,  $\lambda_1(\sT)$ is a simple eigenvalue. Its ground state can be chosen positive, thanks to Sturm-Liouville theory, and we denote the positive and normalized ground state by $\psi$. Next, we provide a full asymptotic expansion of $\lambda_1(\sT)$,  and an approximation of $\psi$.
\begin{proposition}\label{prop:harm-app*}
Suppose that $v$ satisfies \eqref{eq:cond-v}. Then, there exists a sequence of real numbers $(\mu_k)_{k\geq0}$ with $\mu_0=v(0)$ and $\mu_1=\sqrt{2v''(0)}\bigl(1+e_0 \bigr)$,  such that the ground state energy of $\sT$ satisfies, for every integer $N\geq 2$, as $h$ tends to $0$, 
\[\lambda_1(\sT)
\underset{h\to0}
{=} \sum_{k=0}^N\mu_kh^k+\mathcal O(h^{N+1}).\]

Moreover, for every $C_2>C_1>0$, the positive normalized ground state $\psi$ 
  satisfies
\[ 
\begin{gathered}
\|\psi - \psi^\harm\|_{L^\infty(C_1\sqrt{h},C_2\sqrt{h})}=\cO(h^{1/2} ),\\
\|(\psi - \psi^\harm)' \|_{L^\infty(C_1\sqrt{h},C_2\sqrt{h})}=\cO(1 ),
\end{gathered} \]
where $\psi^\harm$ is the normalized function defined as
\begin{equation}\label{eq:qm-h-app}
\psi^\harm(r) =h^{-1/2} u_1(r/\sqrt{h}),
\end{equation}
and $u_1=u^\harm$ is the function introduced in \eqref{eq:def-u1}, with $\beta=v''(0)/2$.
\end{proposition}
\begin{proof}
Fix an integer $N\geq 2$ and expand the potential $v$ by a Taylor series 
\[v(r)=v(0)+\frac{v''(0)}{2!}r^2+\cdots+\frac{v^{(2N)}(0)}{(2N)!}r^{2N}+\cO(r^{2N+2}).\] Applying the scaling $r\mapsto h^{1/2}r$, we get that $\sT$  is unitarily equivalent to
\[\begin{aligned}
\widehat \sT&:=-\frac{\dd^2}{\,\dd r^2}-\frac1r\frac{\dd}{\,\dd r}+\frac{e_0^2}{r^2}+h^{-1}\bigl(v(h^{1/2}r)-v(0)\bigr)\\
&=\widehat \sT_1+h\widehat \sT_2+\cdots+h^N\widehat \sT_N+h^{N+1}R_N,
\end{aligned}\]
where
\[\widehat \sT_1=-\frac{\dd^2}{\,\dd r^2}-\frac1r\frac{\dd}{\,\dd r}+\frac{v''(0)}2 r^2+\frac{e_0^2}{r^2},\]
and $\widehat \sT_2,\cdots,\widehat \sT_N,R_N$ are multiplication operators,
\[\widehat \sT_k=\frac{v^{(2k)}(0)}{(2k)!}r^{2k},\quad R_N=\cO(r^{2N+2}).\]
With $\beta=v''(0)/2$, 
the operator $\widehat \sT_1$ is $\sT^\harm$ studied earlier. We choose $
\mu_1=2\sqrt{\beta}\bigl(1+e_0\bigr)$ the ground state energy of $\widehat \sT_1$, and we choose $
u_1=u^\harm$ its corresponding normalized ground state,
 as in \eqref{eq:def-u1}.

 The operator $\widehat \sT_1-\mu_1$ being invertible on the orthogonal complement of $u_1$,  we can choose real numbers $(\mu_k)_{k=2}^N$ and functions $(u_k)_{k=2}^N$ in the domain  of $\widehat \sT_1$ satisfying
\[ \begin{gathered}
(\widehat \sT_1-\mu_1)u_2+(\widehat \sT_2-\mu_2)u_1=0\\
(\widehat \sT_1-\mu_1)u_3+(\widehat \sT_2-\mu_2)u_2+(\hat T_3-\mu_3)u_1=0\\
(\widehat \sT_1-\mu_1)u_4+(\widehat \sT_2-\mu_2)u_3+(\widehat \sT_3-\mu_3)u_2+(\widehat \sT_4-\mu_4)u_1=0\\
\cdots
\end{gathered}\]
In fact,  we take
\[\mu_2:=\langle \widehat \sT_2u_1,u_1\rangle=\frac{(e_0+1)(e_0+2)v^{(4)}(0)}{4!\, \beta},\]
so that $(\widehat \sT_2-\mu_2)u_1$ is orthogonal to $u_1$  in $L^2(\mathbb R^+; rdr)$. Then, we choose the function $u_2$ orthogonal to $u_1$ such that 
\[ (\widehat \sT_1-\mu_1)u_2=-(\widehat \sT_2-\mu_2)u_1.\]
Writing $u_2=b_2u_1$, the function $b_2$ should satisfy
\[b_2''(r)+\frac1{r}b_2'(r)=\widehat \sT_2 -\mu_2=\frac{v^{(4)}(0)}{4!}r^4-\mu_2,\]
and the general solution of this equation is
\[b_2(r)=\frac{v^{(4)}(0)}{4!}\frac{r^6}{30}-\mu_2\frac{r^2}{4}+c_1\ln r+c_2, \]
with constants $c_1,c_2$. For $u_2$ to be in the domain of $\hat \sT_1$, we choose $c_1=0$ and get
\[u_2(r)=\Bigl(\frac{v^{(4)}(0)}{4!}\frac{r^6}{30}-\mu_2\frac{r^2}{4}+c_2\Bigr)u_1.\]
Then, we choose $c_2$ so that $u_2$ is orthogonal to $u_1$ in $L^2(\R_+,r\,\dd r)$.
The construction of $\mu_3,u_3,\cdots$ follows the same process.

With $u=u_1+hu_2+\cdots+h^Nu_N$ and $\mu=\mu_1+\mu_2h+\cdots+\mu_{N}h^{N-1}$,  we have
\[\|(\widehat \sT-\mu)u\|_{L^2(\R_+,r\,\dd r)}=\cO(h^{N}),\quad \|u\|_{L^2(\R_+,r\,\dd r)}= \|u_1\|_{L^2(\R_+,r\,\dd r)}+\cO(h),\]
and by the spectral  theorem, the ground state energy $\lambda_1(\widehat \sT)$ of $\widehat \sT$ satisfies $\lambda_1(\hat \sT)=\mu+\cO(h^N)$. Therefore, we have 
\[\lambda_1(\sT)-v(0)=h\lambda_1(\widehat \sT)= h\mu+\cO(h^{N+1}).\]
Finally, denoting by $\hat\psi$ the positive and normalized ground state of $\widehat \sT$, and by $\hat u$ the orthogonal projection of $u$ on $\hat\psi$, we have by Proposition~\ref{prop:harm-app} that $\lambda_2(\widehat \sT)-\lambda_1(\widehat \sT)\geq 2\sqrt{\beta}$,  hence by a routine application of the spectral theorem, 
\[ \|u-\hat u\|_{L^2(\R_+,r\,\dd r)}=\cO(h^N),\]
and consequently, since  $0 < C_1 < C_2<+\infty$, 
\[ \int_{\frac{1}2C_1}^{2C_2} \bigl(|(u-\hat u)''|^2+|(u-\hat u)'|^2+|u-\hat u|^2 \bigr)\, \,\dd r=\cO(h^{2N}).\]
 
By Sobolev embedding and a change of variable, we get the result concerning the approximation of $\psi$ in $[C_1\sqrt{h},C_2\sqrt{h}]$.
\end{proof}
\begin{remark}\label{rem:harm-wkb} Let 
\begin{subequations} \label{eq:dd0}
\begin{equation}
d_0(r)=\sqrt{2v''(0)}\, r^2/4\,.
\end{equation}
Since  the even extension of $v$ is smooth, the Agmon distance in  \eqref{eq:dist} can be extended as an even $C^\infty$-function on $\mathbb R$, and using $v'(0)=v^{(3)}(0)=0$, it satisfies as  $r\to0$,
\begin{equation} d(r)=d_0(r)+\mathcal O(r^4),\quad d'(r)=d_0'(r)+\mathcal O(r^3),\quad d''(r)=d_0''(r)+\mathcal O(r^2).
\end{equation}
\end{subequations}
Consequently,  in the interval $[0,C\sqrt{h}]$ and with  $\beta=\sqrt{v''(0)/2}$, we can express $\psi^\harm$ and its derivative  in terms of the Agmon distance  as
\[\begin{gathered}
\psi^\harm(r)=h^{-\frac{1+e_0}{2}} \bigl(A_0+\cO(h) \bigr) r^{e_0}\ee^{-d(r)/h}, \\
(\psi^\harm)'(r)=h^{-\frac{1+e_0}{2}} \Bigl(A_0\frac{\dd}{\,\dd r}\Bigl( r^{e_0}\ee^{-d(r)/h}\Bigr) +\cO(h^{1/2})\Bigr),
\end{gathered}\]
where $A_0$ is the constant introduced in \eqref{eq:def-hat-a}.
\end{remark}

\subsection{Decay at infinity}

It is standard \cite{He88} to derive a rough decay estimate of the ground state $\psi$ as in the following:

\begin{proposition}\label{prop:dec-psi}
Suppose that $v$ satisfies \eqref{eq:cond-v}.  Then, for any $\delta\in(0,1)$, the ground state energy $\psi$ of $\sT$ satisfies as $h\to0$, 
\[\int_{\R_+}\Bigl( h|(\ee^{(1-\delta)d(r)/h}\psi)'|^2+ \Bigl(\frac{he_0^2}{r^2}+1\Bigr)|\ee^{(1-\delta)d(r)/h}\psi(r)|^2\Bigr)r\,\dd r=\cO(1).\]
 where $d$ is as in \eqref{eq:dist}. In particular, for any $R>0$, we have,
\[ \|\ee^{(1-\delta)d(r)/h}\psi \|_{H^1(R;+\infty)}=\cO(h^{-1/2}). \]
\end{proposition}
\begin{proof}
We denote by $\langle\cdot,\cdot\rangle$ and $\|\cdot\|$ the inner product and norm in $L^2(\R_+,r\,\dd r)$. Let $\delta\in(0,1)$ and $n\in\N$. Define the function $\zeta$ on $\R_+$ as 
\[\zeta(r)=(1-\delta)\min(d(r),n).\] Then, with $\sq$ defined in \eqref{eq:def-form-T0},   we have the  following identity
\[\begin{aligned}\sq\bigl(\ee^{\zeta/h}\psi\bigr)-\bigl\| \zeta'\ee^{\zeta/h}\psi\bigr\|^2&=\bigl\langle \sT\psi,\ee^{2(1-\delta)\zeta/h}\psi\bigr\rangle\\
&=\lambda_1(\sT)\bigl\| \ee^{\zeta/h}\psi\bigr\|^2.
\end{aligned}\]
Noticing that $|\zeta'(r)|\leq (1-\delta)d'(r)$, and using Proposition~\ref{prop:harm-app}, we get 
\[\int_{\R_+}\Bigl(h^2|(\ee^{\zeta/h}\psi)'|^2+\Bigl(\frac{h^2e_0^2}{r^2}+\delta(v(r)-v(0))+\cO(h)\Bigr)|\ee^{\zeta/h}\psi|^2\Bigr)r\,\dd r\leq 0. \] 
Choose $C_0,M_0,h_0>0$ such that, for $d(r)\geq C_0h$  and $h\in(0,h_0]$, we have 
\[\delta(v(r)-v(0))+\cO(h)\geq M_0 h\,.\] Consequently, we have, for $h\in(0,h_0]$ and a positive constant $\widehat M_0$,
\[ h^2\int_{\R_+}\Bigl(|(\ee^{\zeta/h}\psi)'|^2+\frac{e_0^2}{r^2}|\ee^{\zeta/h}\psi|^2\Bigr) r\,\dd r+M_0h\int_{C_0h \leq d(r)\leq n} |\ee^{\zeta/h}\psi|^2 r\,\dd r\leq \widehat M_0h. \]
We rewrite this estimate as
\[ \int_{\{d(r)\leq n\}}\Bigl(h^2|(\ee^{(1-\delta)d(r)/h}\psi)'|^2+\Bigl(\frac{h^2e_0^2}{r^2}+M_0h\Bigr)|\ee^{(1-\delta)d(r)/h}\psi|^2\Bigr)r\,\dd r\leq \widetilde M_0h. \]
The constants $\delta,M_0,\widetilde M_0,h_0$ being independent of $n$, we get by monotone convergence, 
\[ \int_{\R_+}\Bigl(h^2|(\ee^{(1-\delta)d(r)/h}\psi)'|^2+\Bigl(\frac{h^2e_0^2}{r^2}+M_0h\Bigr)|\ee^{(1-\delta)d(r)/h}\psi|^2\Bigr)r\,\dd r\leq \widetilde  M_0h. \]
\end{proof}

\subsection{WKB approximation}

\subsubsection{WKB construction}

We recall from \eqref{eq:dist} that 
\begin{equation*}
d(r)=\int_0^r\sqrt{v(\rho)-v(0)}\,\dd \rho.
\end{equation*}

With the help of $d$, we have a useful representation of the operator $\sT$ acting on radial functions:
\begin{lemma}\label{lem:op-T0} 
Let $\mu_0=v(0)$ and $d$ be as above. If $g\in C^\infty(\R_+)$ is a radial function, then
\[
\ee^{d(r)/h}(\sT-\mu_0)\ee^{-d(r)/h}g=(h\cL_1 +h^2\cL_2) g,
\]
where
\[\cL_1:=2d'(r)\frac{\dd}{\,\dd r}+d''(r)+\frac{d'(r)}{r},\quad \cL_2:=-\frac{\dd^2}{\,\dd r^2}-\frac1r\frac{\dd}{\,\dd r}+\frac{e_0^2}{r^2}.
\]
Moreover,
\[\begin{gathered}
 \widehat\cL_1:=r^{-e_0}\cL_1r^{e_0}=2d'(r)\frac{\dd}{\,\dd r}+(1+2e_0)\frac{d'(r)}{r}+d''(r),\\
 \widehat\cL_2:=r^{-e_0}\cL_2r^{e_0}=-\frac{\dd^2}{\,\dd r^2}-\frac{1+2e_0}{r}\frac{\dd}{\,\dd r}.
 \end{gathered}\]
\end{lemma}

We will consider a special class of radial functions in the domain of the operator $\sT$. 
Recall that a  ($h$ independent) function $\widehat a$ on $\R_+$ is said to be  at most  exponentially growing if there is $m>0$ such that $\widehat a(r)=\cO(\ee^{mr})$ in a neighborhood of $+\infty$.
\begin{lemma}\label{lem:dom-T0}
Suppose that $\widehat a\in C^2(\overline\R_+)$ and $\widehat a,\widehat a',\widehat a'' $ are at most exponentially growing. Then, there exists $h_0>0$ such that, for all $h\in(0,h_0]$, the function defined as
\[\psi(r)=r^{e_0}\widehat a(r)\ee^{-d(r)/h}\]
belongs to the domain of $\sT$.
\end{lemma}
\begin{proof}
By Lemma~\ref{lem:op-T0}, we have
\[\sT \psi= \ee^{-d(r)/h} r^{e_0}\bigl(h\widehat\cL_1+h^2\widehat\cL_2 \bigr)\widehat a. \] 
Since $\widehat a\in C^2(\overline\R_+)$ we know that 
$\widehat a' $ and $\widehat a''$ are bounded in a neighborhood of $r=0$. Moreover, $d'(0)=0$ and $d'(r)/r$ is bounded in a neighborhood of $0$. Thanks to the smoothness and exponential growth assumptions on $\widehat a$, the function   
$\bigl(h\widehat\cL_1+h^2\widehat\cL_2 \bigr)\widehat a$ grows exponentially. In a neighborhood of $+\infty$, the function $d(r)$ satisfies
\[
d(r)=d(a)+|v(0)|^{1/2}r\,,\]
hence, for $h$ sufficiently large,  $\sT\psi\in L^2(\R_+)$. In a similar fashion, we have that $\psi,\psi/r,\psi'\in L^2(\R_+)$. This prove that $\psi\in\mathrm{Dom}(\sT)$.
\end{proof}
\begin{remark}\label{rem:dom-T0}
The function $u=(r^{e_0}\ln r)\ee^{-d(r)/h}$ does not belong to the domain of $\sT$,  due to the singularity at $0$.
\end{remark}
To perform the WKB construction, we will need to solve a first order differential equation of Fuchs type, where the general solution is explicit.
\begin{lemma}\label{lem:DE-sol}
Let $\mu_1=(1+e_0)\sqrt{2v''(0)}$. Suppose that $g$ is a smooth function on $\R_+$ and consider    the differential equation
\begin{equation}\label{eq:F}\tag{F} (\widehat\cL_1-\mu_1)\widehat a=g, 
\end{equation}
where $\widehat\cL_1$ is introduced in Lemma~\ref{lem:op-T0}.
\begin{enumerate}[\rm i)]
\item If  $g\in C(\overline\R_+)$ and $g(0)=0$, then, for every $C\in\R$, \eqref{eq:F} has a unique smooth solution on $\overline\R_+$ such that $\widehat a(0)=C$, and if $g=0$, it holds
\[\widehat a(r)=C \exp\Bigl(-\int_0^rp_{e_0}(\rho)\dd\rho \Bigr),\]
where $p_{e_0}$ is introduced in \eqref{eq:def-phase}.
\item  If $g(0)\not=0$,  \eqref{eq:F} has a smooth solution $\widehat a$ in $\R_+$ such that if $g(r)=r^{-k}$ with $k\geq 1$, then $\widehat a$ satisfies
\[r^k\widehat a(r)=\cO(1)\quad (r\to0).  \]
\item 
 Any solution of \eqref{eq:F} satisfies, for any $r\geq R_0\geq \sigma$, 
\[\widehat a(r)=\Bigl(\frac{R_0}{r}\Bigr)^{e_0+\frac12}\ee^{\zeta_0(r-R_0)}\widehat a(R_0)+ \frac{r^{-e_0-\frac12}\ee^{\zeta_0r}}{2\sqrt{|v(0)|}}\int_{R_0}^rg(\rho)\rho^{e_0+\frac12}\ee^{-\zeta_0\rho}\dd\rho, \]
with $\zeta_0=(1+e_0)\sqrt{\frac{v''(0)}{2|v(0)|}}$. \end{enumerate}
\end{lemma}
\begin{proof}
We have (see \eqref{eq:dd0}),
\[\begin{gathered}
2d'(r)= \sqrt{2v''(0)}r+\cO(r^2),\\
b(r):=(1+2e_0)\frac{d'(r)}{r}+d''(r)= (1+e_0)\sqrt{2v''(0)}+\cO(r^2).
\end{gathered}\]
We can prove that the functions 
\[p(r)=\frac{b(r)-\mu_1}{2d'(r)}\quad\mbox{and}\quad \frac{p(r)}{r}\]
 extend to  smooth functions  on $\overline\R_+$, and consequently we can introduce the smooth function
\[m(r)=\int_0^r p(\rho)\dd\rho.\]
Assuming that $g(0)=0$, the function $q(r)=\frac{g(r)}{2d'(r)}$ can be extended to a  smooth function  on $\overline\R_+$ and the differential  equation \eqref{eq:F}
reads as \[\frac{\dd \widehat a}{\,\dd r}+p(r)\widehat a=q(r)\,.
\]
Its solution satisfying $\widehat a(0)=C$ is 
\[\widehat a(r)=\left( C+\int_0^rq(\rho)\ee^{m(\rho)}\dd\rho\right)\ee^{-m(r)}.\]

If $g(0)\not=0$, a smooth solution of \eqref{eq:F} is
\[\widehat a(r)= \left(-\int_r^1q(\rho)\ee^{m(\rho)}\dd\rho\right)\ee^{-m(r)},\]
and if $g(r)=r^{-k}$ with $k\geq 1$, we have $r^k\widehat a(r)=\cO(1)$ as $r\to0$.\medskip

In $[R_0,+\infty)$, since $v$ vanishes, we have
\[\begin{aligned}
 m(r)&=m(R_0)+\int_{R_0}^{r}\Bigl(\frac{1+2e_0}{2\rho}-(1+e_0)\sqrt{\frac{v''(0)}{2|v(0)|}}\,\Bigr)\dd \rho\\
 &=m(R_0)+(e_0+1/2)(\ln r -\ln R_0)-(1+e_0)\sqrt{\frac{v''(0)}{2|v(0)|}}(r-R_0),
 \end{aligned}\] 
Consequently, with $\zeta_0=(1+e_0)\sqrt{\frac{v''(0)}{2|v(0)|}}$, we have
\[\ee^{m(r)}=\Bigl(\frac{r}{R_0}\Bigr)^{e_0+\frac12}\ee^{-\zeta_0(r-R_0)}\ee^{m(R_0)}\mbox{ for }r\geq R_0. \]
Solving \eqref{eq:F} in the interval $[R_0,+\infty)$  yields the last formula in Lemma~\ref{lem:DE-sol}.
\end{proof}
We are now  ready to perform the WKB construction.  Notice that the construction for $e_0>0$ is technically different from $e_0=0$, since  the  WKB ansatz is in general singular at $r=0$, and the construction provides us just with the terms $\mu_0$ and $\mu_1$, while the remaining terms of the sequence are taken as in Proposition~\ref{prop:harm-app*}. This technical difference is behind the difficulty to deal with the Aharonov-Bohm potential, which introduces a singular term in the differential equations we encounter. A similar difficulty  (the non existence of a WKB solution on $\overline{\mathbb R}_+$)  was encountered in the context of degenerate potential wells by Martinez-Rouleux \cite{MR}.
\begin{proposition}\label{prop:WKB-sol}
Let  $(\mu_k)_{k\geq 0}$ be the sequence as in Proposition~\ref{prop:harm-app*}. There exists a sequence $(\widehat a_k)_{k\geq 0}$ in $C^\infty(\R_+)$ such that the following holds:
\begin{enumerate}[a)]
\item The function $\widehat a_0$ is given  in \eqref{eq:def-hat-a}, and $\widehat a_0,\widehat a_1\in C^\infty(\overline\R_+)$. 
\item For every $k\geq 0$, the function $\widehat a_k$ and its derivatives grow polynomially, and for $k\geq 2$,  $\widehat a_k(r)=\cO(1/r^{k-1})$ as $r\to0$.
\item\label{est-b)} For every integer $N\geq0$, the function defined on $\R_+$ as
\[\psi_N^\wkb(r):=h^{-\frac{1+e_0}{2}}r^{e_0}\Bigl(\sum_{k=0}^N h^k\widehat a_k(r)\Bigr)\ee^{-d(r)/h} ,\]
satisfies  uniformly on every bounded interval of $\R_+$, 
\[
w_N(r) \Biggl|\ee^{d(r)/h}\Biggl(\sT-\sum_{k=0}^{N+1}\mu_k h^k\Biggr)\psi_N^{\wkb}\Biggr|=\mathcal O\bigl(h^{N+\frac{3-e_0}{2}}r^{e_0}\bigr),
\]
 with $w_0(r)=w_1(r)=1$, and $w_N(r)=r^{N+1}$ for $N\geq 2$.
\item The functions $\psi_0^\wkb,\psi_1^\wkb$ belong to the domain of the operator $\sT$, and we have
\[
\|\psi_0^\wkb\|=1+o(1).\]
\item For $e_0=0$, $\psi_N^\wkb$ belongs to the domain of $\sT$ for every integer $N$, and we have
 uniformly on every bounded interval of  $\R_+$,
\[
\Biggl|\ee^{d(r)/h}\Biggl(\sT-\sum_{k=0}^{N+1}\mu_k h^k\Biggr)\psi_N^{\wkb}\Biggr|=\mathcal O\bigl(h^{N+\frac{3}{2}}\bigr).
\]
\end{enumerate}
\end{proposition}
\begin{proof}
Consider a sequence of functions $ (\widehat a_k)_{k\geq 0}\subset C^\infty(\R_+)$.
For all $k\geq 0$ and $r>0$, put $a_k(r)=r^{e_0}\widehat a_k(r)$.

For any  $N \in \mathbb N $,   Lemma~\ref{lem:op-T0} yields
\[
\ee^{d(r)/h}\Bigl(\sT_0-\sum_{k=0}^{N+1}\mu_kh^k\Bigr)\ee^{-d(r)/h}\sum_{k=0}^N h^ka_k=
r^{e_0}\cE_N(r),
\]
where
\[
\begin{aligned}
\cE_N&=\Bigl(h\widehat\cL_1+h^2\widehat\cL_2-\sum_{k=1}^{N+1}\mu_kh^k\Big)\sum_{k=0}^N h^k\widehat a_k\\
&=h(\widehat\cL_1-\mu_1)\widehat a_0\\
&\quad+h\sum_{k=1}^{N} h^{k}\Bigl((\widehat\cL_1-\mu_1)\widehat a_k+\widehat\cL_2\widehat a_{k-1}-\sum_{\ell =2}^{k+1}\mu_\ell \widehat a_{k+1-\ell}\Bigr)+h^{N+2}\cR_N,
\end{aligned}
\]
and
\[\cR_N=\begin{cases}
(\widehat\cL_2-\mu_2)\widehat a_{N} -\sum\limits_{k=3}^{N+1}\mu_k\widehat a_{N-k+2}-h\mu_{N+1}\widehat a_N&\mbox{if }N\geq 2\\
(\widehat\cL_2-\mu_2)\widehat a_{1}&\mbox{if }N=1.
\end{cases}\]

\paragraph{\bf Definition of $\widehat a_0$.}~\medskip

To get the cancellation of  the term of order $h$, we choose $\widehat a_0$ such that \[(\widehat\cL_1-\mu_1)\widehat a_0=0\,,\]
 which is possible by Lemma~\ref{lem:DE-sol}. In fact, we have 
\[\widehat a_0(r)=C \exp\left( -\int_0^r p(\rho)\dd\rho\right), \]
where 
\[p(r)=\frac{d''(r)+(1+2e_0)\frac{d'(r)}{r}-\sqrt{2v''(0)}(1+e_0)}{2d'(r)},\]
and $C>0$ is a constant to be determined below. \\ Notice that \[\widehat a_0(0)=C \mbox{ and } \widehat a'_0(0)=0\,.\]
To determine  the constant $C$ appearing  in the definition of $\widehat a_0$,
we compute
\[
\begin{aligned}
\int_{\R_+} |\psi_0^\wkb(r)|^2r\,\dd r&= h^{-1-e_0}\int_0^{+\infty}|\widehat a_0(r)|^2 \ee^{-2d(r)/h} r^{1+2e_0}\,\dd r\\
&\sim  |\widehat a_0(0)|^2\frac{\Gamma(e_0+1)}{2}\left(\frac{2}{v''(0)}\right)^{(1+e_0)/2},\end{aligned}\] 
 and we recall that  the notation $f(h)\sim  g(h)$ means $f(h)=(1+o(1))g(h)$ as $h\to0$. 

In order for the norm of $\psi_0^\wkb$ in  $L^2(\R_+,r\,\dd r)$ to be asymptotically $1$, we choose
\begin{equation} \label{eq:formC} C=\sqrt{\frac{2}{\Gamma(e_0+1)}}\left(\frac{v''(0)}{2}\right)^{\frac{1+e_0}4}  .\end{equation}
\paragraph{\bf Definition of $\widehat a_1$.}~\medskip

To get the cancellation of  the term of order $h^2$, we choose $\widehat a_1$  so that
\[ (\widehat\cL_1-\mu_1)\widehat a_1+(\widehat \cL_2-\widehat\mu_2)\widehat a_0=0.\]
This is possible  if we choose  $\widehat\mu_2$ as
\[\widehat\mu_2= \widehat\cL_2\widehat a_0|_{r=0},\]
and apply Lemma~\ref{lem:DE-sol} with 
\[g=g_1:=-(\widehat \cL_2-\widehat\mu_2)\widehat a_0=\Bigl( p(r)^2-p'(r)+\frac{(1+2e_0)p(r)}{r}+\widehat\mu_2\Bigr)\widehat a_0.\] We observe that $g_1$ is smooth and we get a smooth solution $\widehat a_1$ with $\widehat a_1(0)=0$ and such that $\widehat a_1$ and its derivatives grow exponentially. Moreover,  we have the explicit expressions
\[ \widehat a_1=w_1\widehat a_0,\quad w_1(r)=C^{-1}\int_0^r\frac{g_1(\rho)}{2d'(\rho)}\ee^{m(\rho)}\dd \rho,\quad m(r)=\int_0^r p(\rho)\dd\rho.\]
Thanks to Lemma~\ref{lem:dom-T0}, we know that that the functions $\phi_0^\wkb,\phi_1^\wkb$ belong to the domain of the operator $\sT$. 
\begin{remark}
\begin{enumerate}[i)]~\medskip
\item Since $\psi_1^\wkb$ is in the domain of $\sT$,  we get from the above construction and the spectral theorem that
\[\lambda_1(\sT)=v(0)+\sqrt{2v''(0)}(1+e_0)h+\widehat\mu_2 h^2+\cO(h^3).\] 
Proposition~\ref{prop:harm-app*} then implies $\widehat\mu_2=\mu_2$.
\item In the case $e_0=0$, we have $\widehat\mu_2=0$ and $g_1=\cO(r^2)$.
\end{enumerate} 
\end{remark}

So far, we have constructed $\widehat a_0,\widehat a_1$, with $\cE_1=h^3\cR_1$ and $\cR_1=\cO(1)$ uniformly on every bounded interval of $\R_+$. \medskip

\paragraph{\bf Definition of $\widehat a_2$.}~\medskip

To  cancel the term of order $h^3$, we choose $\widehat a_2$ such that
\[(\widehat\cL_1-\mu_1)\widehat a_2+(\widehat\cL_2-\mu_2)\widehat a_1-\mu_3\widehat a_0=0,\]
and recall that $\mu_1,\mu_2,\mu_3$ are as in Proposition~\ref{prop:harm-app*}.
Let
\[g_2=\mu_3\widehat a_0-(\widehat\cL_2-\mu_2)\widehat a_1,\]
and notice that
\begin{multline*}
 -(\widehat\cL_2-\mu_2)\widehat a_1=\Bigl(p(r)^2-p'(r)+\frac{(1+2e_0)p(r)}{r} +\mu_2\Bigr)\widehat a_1\\
 +q_1'(r)-p(r)q_1(r)+\frac{(1+2e_0)q_1(r)}{r},
 \end{multline*}
 where
 \[q_1(r)=\frac{g_1(r)}{2d'(r)}.\]
 The function $q_1$ is smooth, but $q_1(0)\not=0$, so 
\begin{multline*}
g_2(r)=\mu_3\widehat a_0(r)+q_1'(r)-p(r)q_1(r)\\+\Bigl(p(r)^2-p'(r)+\frac{(1+2e_0)p(r)}{r} +\mu_2\Bigr)\widehat a_1(r)
+\frac{(1+2e_0)q_1(r)}{r},
\end{multline*}
is singular at $r=0$.  
 We rewrite $g_2$ as
\[g_2(r)=g_2^{\rm reg}(r)+g_2^{\rm sing}(r), \]
where
\[ g_2^{\rm reg}(r)=f_2(r)-C_2,\quad C_2=f_2(0),\quad g_2^{\rm sing}(r)=C_2+\frac{(1+2e_0)q_1(0)}{r},\]
and for $r>0$, $f_2(r)$ is defined as
\begin{multline*}
f_2(r)=\mu_3\widehat a_0(r)+q_1'(r)-p(r)q_1(r)\\
+\Bigl(p(r)^2-p'(r)+\frac{(1+2e_0)p(r)}{r} +\mu_2\Bigr)\widehat a_1(r)+(1+2e_0)\frac{q_1(r)-q_1(0)}{r}.
\end{multline*}
Notice that $g_2^{\rm reg}$ is smooth and satisfies $g_2^{\rm reg}(0)=0$. We then choose $\widehat a_2$ as
\[\widehat a_2=\widehat a_2^{\rm reg}+\widehat a_2^{\rm sing},\]
where $(\widehat\cL_1-\mu_1)\widehat a_2^{\rm reg}=g_2^{\rm reg}$ and  $(\widehat\cL_1-\mu_1)\widehat a_2^{\rm sing}=g_2^{\rm sing}$, and we take $\widehat a_2^{\rm reg}, \widehat a_2^{\rm sing}$ as in Lemma~\ref{lem:DE-sol}.
The function $r^{e_0}\widehat a_2\ee^{-d(r)/h}$ does not belong to the domain of $\sT$ (see Remark~\ref{rem:dom-T0}). Consequently, $\psi_2^\wkb$ does not belong to the domain of $\sT$.
However, we notice that $\cE_2=h^4\cR_2$ and, for every $M>0$, we have 
\[\cR_2=(\widehat\cL_2-\mu_2)\widehat a_2-\mu_3\widehat a_1- h\mu_3 \widehat a_2=\cO(1/r^3)\mbox{  uniformly on }(0,M].\]
 \paragraph{\it Remark.} In the case $e_0=0$, we observe that $q_1(0)=0$ and we can choose $\widehat a_2$ to be smooth on $\overline\R_+$, as stated in Lemma~\ref{lem:DE-sol}. \medskip
\paragraph{\bf The induction process.} ~\medskip

We  repeat this construction and get, for every $k\geq 3$,
\[(\widehat\cL_1-\mu_1)\widehat a_k+\widehat\cL_2\widehat a_{k-1}-\sum_{\ell =2}^{k+1}\mu_l\widehat a_{k+1-l}=0,\]
and $\widehat a_k=\cO(1/r^{k-1})$, $(\widehat\cL_2-\mu_2)\widehat a_k=\cO(1/r^{k+1})$. \\
Consequently, for $N\geq 2$, the function
\[
 \psi_{N}^\wkb(r)=h^{-\frac{1+e_0}{2}}\ee^{-d(r)/h}\Bigl(r^{e_0}\sum_{k=0}^{N} h^k\widehat a_k(r)\Bigr)\]
satisfies  uniformly on every bounded interval of $\R_+$,
\[ r^{N+1}\Biggl|\ee^{d(r)/h}\Biggl(\sT-\sum_{k=0}^{N+1}\mu_k h^k\Biggr)\psi_N^{\wkb}\Biggr|=\mathcal O\bigl(h^{N+\frac{3-e_0}{2}}\bigr).
\]
\end{proof}

\subsubsection{Approximation of ground states}

We now prove Theorem~\ref{thm:T0} following an analysis similar to the one in \cite[Sec.~2]{MR}.
We start by establishing  a refined decay of the ground state $\psi$ away from $r=0$.

\begin{lemma}\label{lem:dec-MR}
Suppose  that $v$ satisfies \eqref{eq:cond-v}.   Then, there exist $M,C_0,h_0>0$ and a positive integer $n_0$  such that,  for all $h\in(0,h_0]$, the ground state $\psi$ of $\sT$ satisfies
\[|\psi(r)|\leq M h^{-n_0/2} r^{-1}\ee^{-d(r)/h}\mbox{ for }r\geq C_0\sqrt{h},\]
where $d$ is the distance introduced in \eqref{eq:dist}.
\end{lemma}
\begin{proof}
Let \[\tilde\lambda=\lambda(\sT)-v(0)\mbox{ and } \tilde d(r)=\int_0^r \sqrt{(v(\rho)-v(0)-\tilde\lambda)_+}\,\dd\rho\,.\]
 Then, by Proposition~\ref{prop:harm-app*}, \[\tilde\lambda=\mu_1h+\cO(h^2)>0\] and we can choose $C_0>0$ such that 
 \[ v(r)-v(0)-\tilde\lambda>0\mbox{ and } \tilde d(r)=d(r)+\cO(h|\ln h|) \mbox{ on }  I_h:=[C_0\sqrt{h},+\infty)\,.
 \]

By integration by parts, we have for all $r>C_0\sqrt{h}$,
\begin{multline*}
\int_{C_0\sqrt{h}}^r \Bigl(h^2|(\ee^{2\tilde d(\rho)/h}\psi)'|^2+\frac{h^2e_0^2}{\rho^2}|\ee^{2\tilde d(\rho)/h}\psi|^2\Bigr)\rho\dd \rho\\
+\int_{C_0\sqrt{h}}^r 
\bigl(v(\rho)-v(0)-|\tilde d'(\rho)|^2\bigr)|\ee^{2\tilde d(\rho)/h}\psi (\rho)|^2\rho\dd \rho\\
-h^2\Bigl(\rho\ee^{2\tilde d(\rho)/h} \psi'(\rho)\psi (\rho)\Bigr)\Big|_{C_0\sqrt{h}}^{r}\\
=\lambda_1(\sT)\int_{C_0\sqrt{h}}^r |\ee^{2\tilde d(\rho)/h}\psi|^2\rho\dd \rho.
\end{multline*}
Observing that $v(\rho)-v(0)-|\tilde d'(\rho)|^2=\lambda_1(\sT)$ on $I_h$, we deduce that
\[  -r\ee^{2\tilde d(r)/h}\psi'(r)\psi (r)\leq  -C_0\sqrt{h}\, \ee^{2\tilde d(C_0\sqrt{h})/h}(\psi'\psi)(C_0\sqrt{h}).\]
By Proposition~\ref{prop:harm-app*}, 
\[\ee^{2\tilde d(C_0\sqrt{h})/h}(\psi'\psi)(C_0\sqrt{h})=\cO\bigl(h^{-3/2} \bigr)\,.
\]
 Consequently we have
\[ -\psi'(r)\psi (r)\leq Ch^{-n_0} r^{-1} \ee^{-2 d(r)/h}.\]
where $C$ is a positive constant and $n_0$ is a positive integer. 

By Proposition~\ref{prop:dec-psi}, we can choose $R>r$ such that \[|\psi(R)|^2=o(\ee^{-2d(r)/h} )\,.\]
 Integrating on $[r,R]$, we get
\[ |\psi (r)|^2\leq 2Ch^{-n_0}\, \int_r^R \rho^{-1} \ee^{-2 d(\rho)/h}\dd\rho+|\psi(R)|^2.\]
On the other hand there exists a positive constant $\gamma$ such that, for every $r\leq \rho\leq R$, we have 
\[d'(\rho)\geq \gamma r \mbox{ and } d(\rho)-d(r)\geq \gamma r(\rho-r)\,.
\]
 Thus, we can bound from above the integral
\[ 
\begin{aligned}
\int_r^R \rho^{-1} \ee^{-2 d(\rho)/h}\dd\rho&\leq r^{-1} \ee^{-2d(r)/h}\int_r^{+\infty} \ee^{-2\gamma r(\rho-r)/h}\dd\rho\\
&=\cO(r^{-2}\ee^{-2d(r)/h}).
\end{aligned}\]
\end{proof}
\begin{proof}[Proof of Theorem~\ref{thm:T0}]~\medskip

The proof relies on a Wronskian argument inspired by  \cite[Sec.~3]{MR}. We choose a positive integer $N$ and, for the ease of notation,  we write \[\psi_N=\psi_N^\wkb\,,\, \mu=\sum_{k=0}^{N+1} \mu_kh^k \mbox{ and } \lambda=\lambda_1(\sT)\,.\]
We also recall from \eqref{eq:def-T0} that
\[
\sT=\sT_{e_0}:=-h^2\frac{\dd^2}{\,\dd r^2}-\frac{h^2}{r}\frac{\dd}{\,\dd r}+v(r) +\frac{h^2e_0^2}{r^2},
\]
It is  straightforward  to check that
\[\begin{aligned}
h^2\bigl (\psi(r\psi_N')-(r\psi')\psi_N\bigr)'&=r\bigl(\psi_N\sT\psi-\psi\sT\psi_N \bigr)\\
&=r(\lambda-\mu)\psi\psi_N-r\psi(\sT-\mu)\psi^N.
\end{aligned} \]
Fix a positive constant $R$. By Proposition~\ref{prop:harm-app*} and Lemma~\ref{lem:dec-MR}, we have on $I_h=[C_0\sqrt{h},R)$,
\[ r(\lambda-\mu)\psi\psi_N=\cO\Bigl(h^{N-\frac{n_0}2}\ee^{-2d(r)/h}\Bigr).\]
Similarly, by Proposition~\ref{prop:WKB-sol} and Lemma~\ref{lem:dec-MR}, we have on $I_h=[C_0\sqrt{h},R)$,
\[r\psi(\sT-\mu)\psi^N=\cO\Bigl(h^{\frac{N}{2}-\frac{n_0}2}\ee^{-2d(r)/h}\Bigr). \]
Therefore, we have
\[\bigl (\psi(r\psi_N')-(r\psi')\psi_N\bigr)' =\cO\Bigl( h^{\frac{N}{2}-\frac{n_0}2-2}\,\ee^{-2d(r)/h}\Bigr)\mbox{ for } C_0\sqrt{h}\leq r\leq R. \]

 Integrating on $[r, R+1]$,   we get\footnote{ With $\tilde R=R+1$, we choose $\delta\in(0,1)$ such that   $(1-\delta)d(\tilde R)>d(R)$, and we integrate $h^2(r\psi')'$ on $[r,+\infty)$ to get $r\psi'(r)=\cO(h^{-2}\ee^{-(1-\delta)d(\tilde R)/h})=\cO(\ee^{-d(r)/h})$ for $\tilde R\leq r\leq 2\tilde R$.} by Proposition~\ref{prop:dec-psi},
\begin{equation}\label{eq:proof-thm-T0*}
 \psi(r)(r\psi_N'(r))-(r\psi' (r))\psi_N(r)=\cO\Bigl( h^{\frac{N}{2}-\frac{n_0}2-2}\, \ee^{-2d(r)/h}\Bigr)\mbox{ for } C_0\sqrt{h}\leq r\leq R. 
 \end{equation}
Notice that there is a constant $c_1>0$ such that, for $C_0\sqrt{h}\leq r\leq R$, we have 
\[\psi_N(r)\geq c_1 h^{1/2}\ee^{-d(r)/h}\,,\]  which follows from  Proposition~\ref{prop:WKB-sol}  and \eqref{eq:def-hat-a}. Therefore, dividing the identity in \eqref{eq:proof-thm-T0*} by $r|\psi_N|^2$ and integrating over $[t\sqrt{h},r]$, with $t>C_0$, we obtain
\begin{equation}\label{eq:proof-thm-T0**}
\frac{\psi(r)}{\psi_N(r)}-\frac{\psi(t\sqrt{h})}{\psi_N(t\sqrt{h})}=\cO\Bigl(h^{\frac{N}{2}-\frac{n_0}2-\frac32}\Bigr) \mbox{ for } C_0\sqrt{h}\leq r\leq R.
\end{equation}
To finish the proof  of the first estimate in Theorem~\ref{thm:T0}, we notice that  Propositions~\ref{prop:harm-app*} and \ref{prop:WKB-sol} yield
\begin{equation}\label{eq:zz1}
\frac{\psi(t\sqrt{h})}{\psi_N(t\sqrt{h})}\underset{h\to0}{=} \frac{u_1(t)+\cO(h^{1/2})}{A_0t^{e_0}\ee^{-\sqrt{v''(0)/2}t^2/4}+\cO(h^{1/2})}=1+\cO(h^{1/2}), 
\end{equation}
and 
\begin{equation}\label{eq:zz2}
\psi_N(r)=\psi_0(r)+\cO(\sqrt{h}\,\ee^{-d(r)/h}) \mbox{ for }t\sqrt{h}\leq r\leq R\,.
\end{equation}
 Next we prove the estimate concerning $\psi'$. We divide \eqref{eq:proof-thm-T0*} by $r\psi_N$ and use that $\psi_N(r)\geq C_1 h^{1/2}\ee^{-d(r)/h}$ to get for $C_0\sqrt{h}\leq r\leq R$,
\begin{align*}
\psi'(r)-\frac{\psi(r)}{\psi_N(r)}\psi_N'(r)=\cO\Bigl(h^{\frac{N}{2}-\frac{n_0}2-1}\ee^{-d(r)/h}\Bigr).
\end{align*}
Thanks to \eqref{eq:proof-thm-T0**}, we have $\frac{\psi(r)}{\psi_N(r)}=1+\cO(h^{1/2})$.
It suffices then to choose $N$ sufficiently large, and notice that, in a compact interval of $\R_+$, we have
\[ \psi_N'(r)=h^{-\frac{1+e_0}{2}}\Bigl( -\frac{d'(r)}{h} r^{e_0}\widehat a_0(r)+\cO(1)\Bigr)\ee^{-d(r)/h}.\]
\end{proof}

\section{The single well operator}\label{sec:sw-op}

In this section, we study  the operator $\sH_{\alpha,0}$ introduced in \eqref{eq:def-H0}, and we prove Theorem~\ref{thm:sw}.

\subsection{Separation of variables}\label{subsec:sep-var}

Using the identification between $L^2(\R^2)$ and $L^2(\R_+,r\,\dd r)\otimes L^2(\mathbb S^1)$ via Fourier modes, we obtain that $\sH_{\alpha,0}$ is unitarily equivalent to the direct sum of the fiber operators in $L^2(\R_+,r\,\dd r)$,
\[-h^2\partial_r^2-\frac{h^2}r\partial_r+h^2\frac{|m-\alpha/h|^2}{r^2}+v(r), \]
parameterized by $m\in\Z$. 

We assume that $e_0\in[0,\frac12]$ is fixed and that $h$ varies in the set $\cJ_\alpha(e_0)$ introduced in \eqref{eq:def-J}. Minimizing over $m\in\Z$, we obtain that the ground state energy $\lambda^\sw(h)$ of $\sH_{\alpha,0}$ is equal to that of the operator $\sT$ introduced in \eqref{eq:def-T0}, and  the  ground states are  given in polar coordinates as
\[\ee^{\ii m\theta}u(r)\]
where $m$ satisfies $|m-\alpha/h|=e_0$ and $u$ is a ground state of $\sT$.

\subsubsection*{Non-half integer flux}

If $e_0\in[0,1/2)$, there is a unique integer $m_*$ satisfying $|m_*-\alpha/h|=e_0$, which can be expressed as
\begin{equation}\label{eq:def-m}
m_*=m_*(\alpha/h):=\begin{cases}
\lfloor \alpha/h\rfloor&\mbox{if }\alpha/h<\lfloor \alpha/h\rfloor+\frac12,\\
\lfloor \alpha/h\rfloor+1&\mbox{if } \alpha/h>\lfloor \alpha/h\rfloor+\frac12.
\end{cases} 
\end{equation} 
Consequently the ground state energy $\lambda^\sw(h)$ is a simple eigenvalue, with the following normalized ground state in $L^2(\R^2)$, defined in polar coordinates as
\begin{equation}\label{eq:def-gs-AB} 
\phi^\sw (r,\theta) =\pi^{-\frac12}\psi(r)\ee^{\ii m_*\theta},
\end{equation}  
where $\psi$ is the positive and normalized (in $L^2(\R_+,r\,\dd r)$)  ground state of $\sT$. The conclusion i) in  
Theorem~\ref{thm:sw} is now a direct consequence of Theorem~\ref{thm:T0} regarding the ground state $\psi$. 
\subsubsection*{Half-integer flux}

Suppose now that $e_0=1/2$. Then, there are two integers $m_*$ and $m_*+1$ such that 
\[ |m_*-\alpha/h|=|m_*+1-\alpha/h|=e_0\]
and we can express $m_*$ as
\[m_*=m_*(\alpha/h):=\frac{\alpha}{h}-\frac12.\]

Consequently, the ground state energy $\lambda^\sw$ has multiplicity $2$, with the following normalized  ground states
\begin{equation}\label{eq:def-gs-AB*} 
\phi^\sw_1  =\pi^{-\frac12}\psi(r)\ee^{\ii m_*\theta},\quad \phi^\sw_2=\pi^{-\frac12}\psi(r) \ee^{\ii (m_*+1)\theta}\,,
\end{equation} 
 where $\psi$ is the positive and normalized ground state of $\sT$. The set $\{\phi^\sw_1,\phi^\sw_2\}$ is orthonormal set, and the conclusion ii) in  
Theorem~\ref{thm:sw} is now  a direct consequence of Theorem~\ref{thm:T0}.

That the ground state energy has mulktiplicity $2$ can be viewed alternatively (as in \cite{HHOO}) by noticing that when $\alpha/h$ is half-integer,  the function $\rho$ defined in polar coordinates by $\rho(r,\theta)=\ee^{2\ii\alpha\theta/h}$ is smooth on $\R^2\setminus\{0\}$, and the anti-linear operator 
\begin{equation}\label{eq:def-op-K}
\sK:L^2(\R^2)\ni u\mapsto \rho \,\overline{u}\in L^2(\R^2)
\end{equation} satisfies
\begin{equation}\label{eq:def-op-K*}
\sH_{\alpha,0}\sK=\sK\sH_{\alpha,0}\quad\mbox{and}\quad \sK^2=\sI. 
\end{equation} 
Thus,  whenever $\phi $ is an eigenfunction  of $\sH_{\alpha,0}$ associated with $\lambda$,  then $\sK\phi$ is also an eigenfunction.

\subsection{Decay of ground states}

For later use, we need an estimate of the decay of the ground states of $\sH_{\alpha,0}$. We have seen earlier that a normalized ground state of $\sH_{\alpha,0}$ is of the form 
\[\phi(r,\theta)=\pi^{-\frac12}\psi(r)\ee^{\ii m\theta}, \]
with $\psi$ the positive normalized ground state of $\sT$. As a consequence of Proposition~\ref{prop:dec-psi},  we can write Agmon estimates giving the decay of $\phi$ at infinity.
\begin{proposition}\label{prop:dec-phi}
Suppose that $\alpha>0$, $e_0\in[0,1/2]$,  and that $v$ satisfies \eqref{eq:cond-v}. Then, for every positive $\delta<1$, there are constants $C,h_0>0$ such that, for $h\in(0,h_0]$, a normalized ground state $\phi$ of $\sH_{\alpha,0}$ satisfies,
\[\sq_\alpha\bigl(\ee^{(1-\delta)d(|\xv |)/h}\phi\bigr)+h\int_{\R^2}\bigl|\ee^{(1-\delta)d(|\xv |)/h}\phi\bigr|^2\dd\xv  \leq Ch,\]
where $d$ is the distance introduced in \eqref{eq:dist}, and $\sq_\alpha$ is the quadratic form introduced in \eqref{eq:def-q-alpha}.
\end{proposition} 

\subsection{Harmonic potential}

In the case where $\vb (\xv )=\beta|\xv |^2$ with $\beta>0$, the spectrum of $\sH_{\alpha,0}$ is purely discrete,  and consists of the discrete eigenvalues  \cite[Thm.~2.1]{H},
\begin{equation}\label{eq:2.2*}
E(m,n)=2h\sqrt{\beta}(1+|m-\alpha/h|+2n),\quad (m\in\Z,~n=0,1,2,\cdots),
\end{equation}
and the multiplicity of a given eigenvalue equals the number of times it is repeated as $m$ and  $n$ vary.
Assuming that $h\in\cJ_\alpha(e_0)$, the lowest eigenvalue is
\begin{equation}\label{eq:gse*}
E_*=2h\sqrt{\beta}\bigl(1+e_0\bigr),
\end{equation}
which is simple when  $e_0\in[0,1/2)$, i.e when $\alpha/h\not\in \Z+\frac12$.  Then,  we can estimate the spectral gap, since  for $m\in\Z$ with $m\not=m_*(\alpha/h)$, we have
\[  |m-\alpha/h|\geq 1-e_0,\]
hence for  $(m,n)\not=(m_*(\alpha/h),0)$, it holds
\begin{equation}\label{eq:gap*}
 E(m,n)-E_*\geq 2h\sqrt{\beta} (1-2e_0)>0. 
 \end{equation}
 If $e_0=1/2$, which amounts to $\alpha/h\in \Z+\frac12$, the lowest eigenvalue has multiplicity $2$, and for $(m,n)\not\in\{((m_*(\alpha/h),0),(m_*(\alpha/h),0)\}$, it holds
  \begin{equation}\label{eq:gap**}
 E(m,n)-E_*\geq 3h\sqrt{\beta} >0. 
 \end{equation}

\subsection{Harmonic approximation}
For later use, we need to approximate the spectral gap for the operator $\sH_{\alpha,0}$ when $v$ satisfies \eqref{eq:cond-v}. Up to small errors, a standard application of  harmonic approximation reduces the problem to the potential $\vb (\xv )=\beta |\xv |^2$, with $\beta=v''(0)/2$. However, we have to distinguish between the case where the lowest eigenvalue of $\sH_{\alpha,0}$ is simple, and the case where it is double.

\subsubsection*{Non half-integer flux}
Suppose that $ e_0\in[0,1/2)$  and that $h\in\cJ_\alpha(e_0)$, where the set $\cJ_\alpha(e_0)$ is as in \eqref{eq:def-J}. We saw earlier that $\lambda^\sw(h)$ is a simple eigenvalue of $\sH_{\alpha,0}$. Let $\lambda_2^\sw(h)$ be the second eigenvalue of $\sH_{\alpha,0}$. Then, it  satisfies
\[\lambda_2^{\sw}(h)-\lambda^{\sw}(h)\geq 
h\sqrt{2v''(0)}\bigl( 1-2e_0\bigr)+h\, \varepsilon(e_0,h)\quad (h\in\mathcal J_\alpha(e_0)), \] 
 where $\varepsilon(e_0,h)\underset{h\to0}{\to} 0$ uniformly with respect to $e_0\in [0,\delta_0)$, for any $\delta_0<\frac12$.\medskip

That the estimate is not uniform with respect to $e_0\in[0,1/2)$ is due to the fact that in the estimate of the remainder, the reciprocal of  $E(m,n)-E_*$ in \eqref{eq:gap*} appears, which is of order $(1-2e_0)^{-1}$.

 \subsubsection*{Half-integer flux}
Suppose now that $ e_0=1/2$  and that $h\in\cJ_\alpha(1/2)$. Then, as explained earlier,  $\lambda^\sw(h)$ is an  eigenvalue of $\sH_{\alpha,0}$ with multiplicity $2$, hence $\lambda_2^\sw(h)=\lambda^\sw(h)$. The third eigenvalue $\lambda_3^\sw(h)$ of $\sH_{\alpha,0}$ $\lambda_3^\sw(h)$ then  satisfies
\[\lambda_3^{\sw}(h)-\lambda^{\sw}(h)\geq 
3h\sqrt{v''(0)/2}+o(h)\quad (h\in\mathcal J_\alpha(1/2)). \] 

\subsection{Application to the double well operator}\label{subsec:dw-op}
We now  consider the double well operator $\sH_\alpha$ introduced in \eqref{eq:def-op}, with the 
potential as in \eqref{eq:def-V}, and we denote by $\{\lambda_n(\sH_\alpha)\}_{n\geq1}$ the  min-max sequence corresponding to $\sH_\alpha$. A standard argument via the min-max principle shows that, up to $o(h)$ errors, the spectrum of $\sH_\alpha$ is given by that of the direct sum $\sH_{\alpha,0}^\ell\oplus\sH_{\alpha,0}^r$, where \[\sH_{\alpha,0}^\ell=(-\ii h\nabla -\alpha\Fb(\xv -\xv_\ell))^2+\vb (\xv -\xv_\ell)\] and
\[\sH_{\alpha,0}^r=(-\ii h\nabla -\alpha\Fb(\xv -\xv_r))^2+\vb (\xv -\xv_r)\,.\]
 The operators $\sH_{\alpha,0}^\ell$ and $\sH_{\alpha,0}^r$ are unitarily equivalent to $\sH_{\alpha,0}$. Consequently, assuming that $h\in \cJ_{\alpha}(e_0)$, we have the following:
\begin{enumerate}[i)]
\item If $e_0\in[0,1/2)$, then the third eigenvalue of $\sH_{\alpha}$ satisfies
\[\lambda_3(\sH_{\alpha})\geq \lambda^{\sw}(h)+
h\sqrt{2v''(0)}\bigl( 1-2e_0\bigr)+o(h), \]
as $h\to0$, uniformly with respect to $e_0\in[0,\delta_0)$ for any $\delta_0<1/2$.
\item   If $e_0=1/2$, then the fifth  eigenvalue of $\sH_{\alpha}$ satisfies 
\[\lambda_5(\sH_\alpha)\geq \lambda^{\sw}(h)+ 
3h\sqrt{v''(0)/2}+o(h). \] 
\end{enumerate}
\section{Quasi-modes for the two wells operator}\label{sec:qm}

In this section, we construct quasi-modes for the operator $\sH_\alpha$ introduced in \eqref{eq:def-op}, which will allow us to reduce the study of the eigenvalue splitting of $\sH_\alpha$ to that of a $2\times2$ matrix for the non half-integer flux (and of a  $4\times4$ matrix for half-integer flux).

In the rest of this section, we suppose that $h$ varies in the set $\cJ_\alpha(e_0)\cap (0,h_0]$ with $e_0\in[0,1/2)$ and $h_0>0$. We return later to the case where $e_0=1/2$, which corresponds to an half-integer flux. 

To lighten the notation, we write  $\sH$ for $\sH_\alpha$, and $\lambda^\sw$ for $\lambda^\sw(h)$.

\subsection{A phase function}

We introduce a function $\eta:\R^2\setminus\{0\}\to[0,2\pi)$ defined via polar coordinates as follows. Every $\xv\in\mathbb R^2\setminus\{0\}$ can be written uniquely in the form $\xv=(r\cos\theta,r\sin\theta)$ with $r=|\xv |>0$ and $\theta\in[0,2\pi)$; we  define $\eta(\xv )=\theta$. 

 Note that $\eta$ is discontinuous on the semi-axis $\Upsilon=\R_+\times\{0\}$, but the  function $\ee^{\ii\eta}$ is smooth on $\R^2\setminus\{0\}$ and  satisfies\footnote{This can be easily checked by using $ \ee^{\ii\eta(\xv)}=x_1/|\xv|+\ii x_2/|\xv|$.}
\[ (-\ii\nabla-\Fb )\ee^{\ii\eta}=0\quad\mbox{on }\R^2\setminus\{0\},\]
where $\Fb$ is the vector field introduced in \eqref{eq:def-F}. Consequently, 
\[\forall u\in C^\infty(\R^2\setminus\{0\}),\quad (-\ii\nabla-\Fb) (\ee^{\ii\eta}u)=\ee^{\ii\eta}(-\ii\nabla u). \]
We will also encounter the function $\ee^{\ii k\eta}$ with $k\in\R$, which is smooth on $\R^2\setminus\{0\}$ if and only if $k\in \Z$;  however, if $k\not \in \Z$, it is smooth on $\R^2\setminus\Upsilon$.

\subsection{Symmetry under reflections}

For later use, we mention first  a symmetry relation under the reflection defined by 
\[\R^2\ni \xv=(x_1,x_2)\mapsto\Phi(\xv )=(-x_1,x_2).\]
 If $\xv\in\R^2\setminus\{0\}$, we have  $\eta(\Phi(\xv ))=\pi-\eta(\xv )$ modulo $2\pi$, hence 
 \[\ee^{-\ii\eta(\xv )}=-\ee^{\ii\eta(\Phi(\xv ))}.\]  
With  $\xv_\ell$ and $\xv_r$  as in \eqref{eq:def-V}, we have $\Phi(\xv )-\xv_\ell=\Phi(\xv -\xv_r)$, and consequently
\begin{equation}\label{eq:symm-eta}
\eta(\xv -\xv_r)=\pi-\eta(\Phi(\xv )-\xv_\ell)\mbox{ modulo }2\pi,\quad\mbox{for }\xv\in\R^2\setminus\{\xv_\ell,\xv_r\}.
\end{equation} 
This gives another  symmetry relation enjoyed by the operator $\sH$. We associate with $\Phi$ the antilinear transformation
\begin{equation}\label{eq:reflection}
L^2(\R^2)\ni u\mapsto \bar\sL_\Phi u\in L^2(\R^2),\quad ( \bar\sL_\Phi u)(\xv )=\overline{u(\Phi(\xv ))}.
\end{equation}
The transformation $\bar\sL_\Phi$  preserves the norm  in $L^2(\R^2)$, but for the inner product we have $\langle\bar\sL_\Phi u, \bar\sL_\Phi v\rangle=\overline{\langle u,v\rangle}$. If $u$ is in the domain of $\sH$, then
\begin{equation}\label{eq:symm-H}
\sH\bar\sL_\Phi u=\bar\sL_\Phi\sH u.
\end{equation} 
 Similarly, we can quantize the symmetry $(x_1,x_2) \mapsto (x_1,-x_2)$ but the most useful is to consider the composition $\Psi$ of the two previous symmetries. 
With 
 $\Psi(\xv)=-\xv$ and \[\sL_\Psi u=u\circ\Psi\,,\]  we have 
 \begin{equation}\label{eq:symm-Ha}
\sH\sL_\Psi =\sL_\Psi \sH\,.
\end{equation} 
\subsection{Construction of a quasi-mode concentrated at  $\xv_\ell$}~\medskip

We introduce the half-plane $\cK_\ell=\{\xv=(x_1,x_2)\in\R^2\colon x_1<L/2\}$ and the function
\begin{equation}\label{eq:ef-ell}
\cK_\ell\ni\xv\mapsto \phi_\ell(\xv )=\phi^{\sw}(\xv -\xv_\ell)\ee^{\ii\alpha\eta(\xv -\xv_r)/h },
\end{equation}
 The phase term $\ee^{\ii\alpha\eta(\xv -\xv_r)/h }$ is smooth on $\cK_\ell$ because for all $\xv \in\cK_\ell$, $\xv-\xv_r\in \R_-\times\R$, hence it does not belong to  $\Upsilon$, the discontinuity set of $\eta$.
Moreover, the phase term amounts to a gauge transformation that allows us to work with the potential $\Fb_\alpha(\xv -\xv_\ell)$ instead of  $\Ab_\alpha$. 
In fact, we  observe that
\begin{equation}\label{eq:ef-ell*}
\sH \phi_\ell=\lambda^{\sw}\phi_\ell+v(\xv -\xv_r)\phi_\ell\quad\mbox{on }\cK_\ell.
\end{equation}
If we write, $x-\xv_\ell=|\xv -\xv_\ell|(\cos\theta,\sin\theta)$, we observe that
\[\phi_\ell(\xv )=\pi^{-\frac12}\psi(|\xv -\xv_\ell|)\exp\Bigg(\ii m_*(\alpha/h)\eta(\xv -\xv_\ell)+\ii\frac{\alpha}{h}\eta(\xv -\xv_r)\Biggr), \]
where $m_*(\alpha/h)$ is the integer  introduced in \eqref{eq:def-m},  and $\psi$ is the positive normalized ground state of the operator $\sT$ introduced in \eqref{eq:def-T0}.
The phase function in \eqref{eq:ef-ell} is not defined in domains containing $\xv_r$, so we  cut it off in a  neighborhood of $\cB_{\sigma}(\xv_r)$,  the disc of center $\xv_r$ and radius $\sigma$, which contains the support of $\vb(\xv-\xv_r)$. 

The cut-off is introduced as follows. Consider a positive number $\varepsilon<\sigma$,  and a  function $\tilde\chi_\ell\in C^\infty(\R;[0,1])$ such that $\tilde\chi_\ell(t)=1$ for $t<L/2-\sigma-\varepsilon$, and $\tilde\chi_\ell=0$ for $t\geq L/2-\sigma$. For $\xv=(x_1,x_2)\in\R^2$, we define 
\[\chi_\ell (\xv )=\tilde\chi_\ell(x_1)\,.\] This yields a function $\chi_\ell\in C^\infty(\R^2;[0,1])$. Finally, put
\begin{equation}\label{eq:qm-ell}
u_\ell(\xv )=\chi_\ell(\xv )\phi_\ell(\xv ).
\end{equation}
\subsection{Quasi-mode concentrated at the right well}
Similarly as we did for the left well, we consider the half-plane 
\[\cK_r=\Phi(\cK_\ell)=\{\xv=(x_1,x_2)\in\R^2\colon x_1>-L/2\}, \]
and the  function defined on $\cK_r$ by
\begin{equation}\label{eq:ef-r}
\phi_r(\xv )=( \bar\sL_\Phi \phi_\ell) (\xv),
\end{equation}
where $ \bar\sL_\Phi$ is the transformation introduced in \eqref{eq:reflection}. 

Moreover, we introduce the cut-off   $\chi_r= \bar\sL_\Phi\chi_\ell\in C^\infty(\R^2;[0,1])$ which  is independent of $x_2$, equals $1$ on $\{x_1>-L/2+\sigma+\varepsilon\}$, and equals $0$ on $\{x_1\leq -L/2 + \sigma\}$. Our 
quasi-mode is then
\begin{equation}\label{eq:qm-r}
u_r(\xv ):=\chi_r(\xv )\phi_r(\xv )= ( \bar\sL_\Phi u_\ell) (\xv ).
\end{equation}

Let us observe that
\begin{equation}\label{eq:qm-p}
|\phi_\ell(\xv )\overline{\phi_r(\xv )}|=\pi^{-1}\psi(|\xv -\xv_\ell|)\psi(|\xv -\xv_r|),
\end{equation}
and on $\{x_1=0\}$,
\[\begin{aligned}
\phi_\ell(\xv)\overline{\phi_r(\xv)}&=\phi_\ell(\xv)^2\\
&=\pi^{-1}\psi(|\xv-\xv_\ell|)^2
\exp\Bigl(2\ii m_*(\alpha/h)\eta(\xv-\xv_\ell)+2\ii\frac{\alpha}{h}\eta(\xv-\xv_r)\Bigr).
\end{aligned} \]
Since $L>0$ and $\eta$ is valued in $[0,2\pi)$, it is straightforward to check that
\[\eta\bigl(-L/2,x_2\bigr))=
\pi-\arctan\Bigl(\frac{x_2}{L/2}\Bigr)\mbox{ for }x_2\in\R, \]
and
\[\eta\bigl(L/2,x_2\bigr)=\begin{cases}
\arctan\bigl(\frac{x_2}{L/2}\bigr)&\mbox{if }x_2>0\\
2\pi+\arctan\bigl(\frac{x_2}{L/2}\bigr)&\mbox{if }x_2\leq 0
\end{cases}.\]
In terms of  
\begin{equation}\label{eq:def-gam0}
\gamma_0=m_*(\alpha/h)-\alpha/h\in\{-e_0,e_0\}\,,
\end{equation}
we have for $\xv=(0,x_2)$,
\[m_*\,\eta(\xv-\xv_\ell)+\frac{\alpha}{h}\eta(\xv-\xv_r)=
\begin{cases}
\pi m_*-\gamma_0\bigl(\pi-\arctan\bigl(\frac{x_2}{L/2}\bigr)\bigr)&\mbox{if }x_2>0\,,\\
3\pi m_*-\gamma_0\bigl(\pi-\arctan\bigl(\frac{x_2}{L/2}\bigr)\bigr)&\mbox{if }x_2\leq0\,.\\
\end{cases} \]
Consequently, with 
\[w(x_2)=\exp\bigl(2\ii\gamma_0\arctan(2x_2/L)\bigr)\,,\] we have, for $\xv=(0,x_2)$,
\begin{equation}\label{eq:qm-p'}
\begin{gathered}
\exp\Bigl(2\ii m_*(\alpha/h)\,\eta(\xv-\xv_\ell)+2\ii\frac{\alpha}{h}\eta(\xv-\xv_r)\Bigr)=
\ee^{-2\ii\pi\gamma_0}w(x_2),\\
\phi_\ell(\xv )\, \overline{\phi_r(\xv )}=\pi^{-1}\ee^{-2\ii\pi\gamma_0}\psi(|\xv -\xv_\ell|)^2\, w(x_2)\, .
\end{gathered}
\end{equation}
\subsection{Some estimates}
It results from Theorem~\ref{thm:T0} and \eqref{eq:qm-p} that, there are positive constants $C$ and $h_0$ such that, 
\begin{equation}\label{eq:qm-p*}
|\phi_\ell(\xv )\overline{\phi_r(\xv )}|\leq C h^{-1-e_0}\ee^{-d(|\xv -\xv_\ell|)/h-d(|\xv -\xv_r|)/h}\quad\mbox{on }\cK_\ell\cap\cK_r.
\end{equation}
Since $\cK_\ell\cap\cK_r=\{\xv=(x_1,x_2)\colon |x_1|<L/2\}\subset\{\xv\colon |\xv-\xv_\ell|\leq L\}$, this leads us to minimize the function
\begin{equation}\label{eq:def-s}
s(\xv ):= d(|\xv -\xv_\ell|)+d(|\xv -\xv_r|)
\end{equation}
for $0\leq |\xv -\xv_\ell|\leq L$. 

Writing $r=|\xv -\xv_\ell|$ and $x-\xv_\ell=r(\cos\theta,\sin\theta)$, we get that
$|\xv -\xv_r|\geq L-r$ and that
\[s(\xv )\geq d(r)+d(L-r)\quad(0\leq r\leq L).\]
Minimizing over $r$, we get
\[s(\xv )\geq 2d(L/2)=2\int_0^{L/2}\sqrt{v(\rho)-v(0)}\dd\rho.\]
Thanks to \eqref{eq:def-s}, we have $s(0)=2d(L/2)$, and we get that
\begin{equation}\label{eq:def-s*}
\min_{|\xv -\xv_\ell|\leq L}s(\xv )=2d(L/2).
\end{equation}
Corning back to \eqref{eq:qm-p*}, we get
\begin{equation}\label{eq:qm-p**}
|\phi_\ell(\xv )\,\overline{\phi_r(\xv )}|\leq C h^{-1-e_0}\ee^{-2d(L/2)/h}\quad\mbox{for }|x_1|< L/2,
\end{equation}
and consequently, we can estimate the inner product in $L^2(\R^2)$ of the quasi-modes $u_\ell$ and $u_r$ introduced in \eqref{eq:qm-ell} and \eqref{eq:qm-r} as follows
\begin{equation}\label{eq:qm-p***}
\langle u_\ell,u_r\rangle=\mathcal O\bigl(h^{-1-e_0}\ee^{-2d(L/2)/h}\bigr).
\end{equation}
Next we estimate the norms, in $L^2(\R^2)$, of the quasi-modes $u_\ell$ and $u_r$. Notice that, due to the symmetry of our construction, we have 
\[\|u_\ell\|^2=\|u_r\|^2=1+\int_{\R^2}(\chi^2_\ell-1)|\phi^\sw(\xv -\xv_\ell)|^2 \dd\xv .\]
 We have \[|\phi^\sw(\xv -\xv_\ell)|=|\psi(|\xv -\xv_\ell|)|\,,\]
 and thanks to Theorem~\ref{thm:T0}, the last integral  is  $\mathcal O\bigl(h^{-1-e_0}\ee^{ -2 d(L-\sigma-\varepsilon)/h} \bigr)$.  Thus,
\begin{equation}\label{eq:qm-norms}
\|u_\ell\|^2=\|u_r\|^2=1+\mathcal O\bigl(h^{-1-e_0}\ee^{-2d(L-\sigma-\varepsilon) /h}\bigr).
\end{equation}
Using \eqref{eq:ef-ell*} and the symmetry relation in \eqref{eq:symm-H}, we get
\begin{equation}\label{eq:qm-operators}
\|(\sH-\lambda^{\sw})u_\ell\|=\|(\sH-\lambda^{\sw})u_r\|=\mathcal O\bigl(h^{-\frac{3+e_0}{2}}\ee^{-d(L-\sigma-\varepsilon)/h}\bigr).
\end{equation}
Finally, we saw earlier that the the third eigenvalue of $\sH_\alpha$ satisfies,
\begin{equation}\label{eq:ev-3rd}
\lambda_3(\sH)\geq \lambda^{\sw}+ch,
\end{equation}
 where $c>0$ is a constant.
\section{Interaction between the wells and the magnitude of tunneling for non half-integer flux}\label{sec:interaction}
Recall that the definition of the quasi-modes $u_\ell$ and $u_r$ involves an arbitrary but fixed constant $\varepsilon\in(0,\sigma)$. Let us fix another positive constant $S_\varepsilon$ such that 
\begin{equation}\label{eq:def-S-ep}
 d(L/2)< S_\varepsilon<\min\bigl(2d(L/2),d(L-\sigma-\varepsilon)\bigr).
\end{equation}
 Let $\cT:=\sH-\lambda^{\sw}$. The construction of the quasi-modes $u_\ell$ and $u_r$ being symmetric, we get
\begin{equation}\label{eq:J1,2}
\langle \cT u_\ell,u_\ell\rangle=\langle \cT u_r,u_r\rangle=:J_0\quad \mbox{and}\quad \langle\cT u_\ell,u_r\rangle=\overline{\langle \cT u_r,u_\ell\rangle}=:J_1.
\end{equation}
Moreover, by \eqref{eq:qm-norms} and \eqref{eq:qm-operators}, we know that 
\[|J_0|=\cO(\ee^{-S_\varepsilon/h} ),\quad |J_1|=\cO(\ee^{-S_\varepsilon/h} ). \]
Thanks to the properties in \eqref{eq:qm-p***}, \eqref{eq:qm-norms}, \eqref{eq:qm-operators} and \eqref{eq:ev-3rd}, there is a basis of the vector space spanned by the eigenfunctions corresponding to the eigenvalues $\lambda_1(\sH)$ and $\lambda_2(\sH)$ such that, the matrix $\sW$  of $\cT$ in this basis can be decomposed as follows (see \cite[Prop.~2.4 \& Corol. 2.5]{HKS})
\[ \sW=\sU+\sR\]
where
\[\sR=\mathcal O(\ee^{-2S_\varepsilon/h})\]
and
\[ \sU=\begin{pmatrix}
\langle \cT u_\ell,u_\ell\rangle&\langle \cT u_\ell,u_r\rangle\\
\langle \cT u_r,u_\ell\rangle&\langle \cT u_r,u_r\rangle
\end{pmatrix}=\begin{pmatrix}
J_0& J_1\\
\overline{J_1}&J_0
\end{pmatrix}.\]
Consequently, computing the eigenvalues of the Hermitian matrix $\sW$, we find a formula for the difference of the first two eigenvalues of $\cT=\sH-\lambda^{\sw}$, and 
\begin{equation}\label{eq:splitting}
\lambda_2(\sH)-\lambda_1(\sH)=2|J_1|+\mathcal O(\ee^{-2S_\varepsilon/h}).
\end{equation}
Next we have to compute the interaction term $J_1$. The key is to express $J_1$ in a pleasant form that will allow us to use \eqref{eq:qm-ell}, thanks to a tricky computation  dating back to \cite{HSj} (see \cite{HSj-pis} and more recently \cite[Sec.~4.2]{FMR} for its formulation in the magnetic context).  We have in view of \eqref{eq:ef-ell*} and  after a partial integration with respect to the first variable,
\[\begin{aligned}
J_1&=\langle [\sH,\chi_\ell]\phi_\ell,\chi_r\phi_r\rangle\\
&=\ii h \int_{\R} \left(P_1\phi_\ell\,\overline{\phi_r}+\phi_\ell\, \overline{P_1\phi_r}\right)(0,x_2)\,\dd  x_2,
\end{aligned}\]
where $P_1=-\ii h\partial_1-A_1$ and $\Ab_\alpha=(A_1,A_2)$ is introduced in \eqref{eq:def-A}.

Notice that
\begin{equation}\label{eq:J1}
J_1=h^2\int_{\R}\bigl(\partial_1\phi_\ell\,\overline{\phi_r}-\phi_\ell\,\overline{\partial_1\phi_r}\bigr)(0,x_2)\, \dd  x_2-2\ii h\int_{\R}(A_1\,\phi_\ell\,\overline{\phi_r})(0,x_2)\,\dd  x_2.
\end{equation}
We will prove the following.
\begin{lemma}\label{lem:J1}
With $\gamma_0$ as in \eqref{eq:def-gam0} and $w(x_2)=\exp\bigl(2\ii\gamma_0\arctan(2x_2/L)\bigr)$, we have
\[J_1=\pi^{-1}h^2\ee^{-2\ii\pi\gamma_0}\int_{\R} w(x_2)\bigl(F(x_2)-\ii G(x_2)\bigr)\dd  x_2,\]
where   $F$ and $G$ are real-valued functions defined by
\begin{equation}\label{eq:def-F,G}
\begin{aligned}
F(x_2)&=\frac{L}{\sqrt{x_2^2+(L/2)^2}}\psi\Bigl(\sqrt{x_2^2+(L/2)^2}\Bigr)\psi'\Bigl(\sqrt{x_2^2+(L/2)^2}\Bigr),\\
G(x_2)&=2\gamma_0\frac{x_2}{x_2^2+(L/2)^2} \Bigl|\psi\Bigl(\sqrt{x_2^2+(L/2)^2}\Bigr)\Bigr|^2.
\end{aligned}
\end{equation}
\end{lemma}
\begin{proof}
In a neighborhood of $\{x_1=0\}$, the functions $\phi_\ell$ and $\phi_r$ can be expressed as follows
\[
\begin{aligned}
\phi_\ell(\xv )&=\pi^{-\frac12}\psi_\ell(\xv )\ee^{\ii m_*(\alpha/h)\eta(\xv -\xv_\ell)}\ee^{\ii\frac{\alpha}{h}\eta(\xv -\xv_r)},\\
\phi_r(\xv )&=\overline{\phi_\ell(\Phi(\xv))}=\overline{\phi_\ell(-x_1,x_2)},
\end{aligned} \]
where
\[\psi_\ell(\xv )=\psi(|\xv -\xv_\ell|),\quad \psi_r(\xv )=\psi(|\xv -\xv_r|). \] 
Consequently,
\[\bigl(\partial_1\phi_\ell\,\overline{\phi_r}-\phi_\ell\,\overline{\partial_1\phi_r}\bigr)(0,x_2)=
2\bigl(\phi_\ell\partial_1\, \phi_\ell\bigr)(0,x_2).\]
To compute $\partial_1\phi_\ell$, we observe that  in $\{|x_1|<L/2\}$, we have (modulo $2\pi$)
\begin{equation}\label{eq:relation-eta}
\begin{aligned}
\eta(\xv -\xv_\ell)&=\arctan\bigl(x_2/(x_1+L/2)\bigr),\\
\eta(\xv -\xv_r)&=\pi+\arctan\bigl(x_2/(x_1-L/2)\bigr),
\end{aligned}
\end{equation}
and consequently
\[
\partial_1\eta(\xv -\xv_\ell)|_{x_1=0}=\partial_1\eta(\xv -\xv_r)|_{x_1=0}=-\frac{x_2}{x_2^2+(L/2)^2}.
\]
Thus, in light of \eqref{eq:qm-p'}, we eventually find
\[
2\pi\ee^{2\ii\gamma_0\pi}\bigl(\phi_\ell\partial_1\phi_\ell\bigr)(0,x_2)
=
w(x_2)F(x_2)
-\ii \gamma_0^{-1}\Bigl(m_*(\alpha/h)+\frac{\alpha}{h} \Bigr)w(x_2)G(x_2).
\]
We still have to handle the second integral in \eqref{eq:J1}. Notice that (see \eqref{eq:def-A}), 
\[A_1(0,x_2)=-\frac{2\alpha x_2}{x_2^2+(L/2)^2},\] 
and,
thanks to \eqref{eq:qm-p'},  we have,
\[
-2\pi \ee^{2\ii\pi\gamma_0} \ii h  A_1\phi_\ell\, \overline{\phi_r}(0,x_2)\\
=2\ii h^2\gamma_0^{-1}\frac{\alpha}{h}\,w(x_2) G(x_2).
\]
It remains to collect the above formulas, insert them into \eqref{eq:J1} and observe that $m_*(\alpha/h)-\alpha/h=\gamma_0$.
\end{proof}
With Lemma~\ref{lem:J1} in hand, we are ready to compute $J_1$.
\begin{lemma}\label{lem:J1*}
With $\gamma_0$ as in \eqref{eq:def-gam0} and $J_1$ as in \eqref{eq:J1,2}, we have
\[ J_1=-\ee^{-2\ii\gamma_0\pi}|J_1|,\]
and
\begin{equation}\label{eq:comp-J1}
|J_1|=\frac12C(L,v,e_0) h^{\frac12-e_0} \ee^{-2d(L/2)/h}+o\Bigr(h^{\frac12-e_0}\ee^{-2d(L/2)/h} \Bigr),
\end{equation}
where $d$ is the distance introduced in \eqref{eq:dist}, $C(L,v,e_0)$ is defined as
\begin{equation}\label{eq:def-C}
C(L,v,e_0)=\frac{4(L/2)^{2e_0+\frac12}|\widehat a_0(L/2)|^2}{\pi^{\frac12}|v(0)|^{\frac14}},
\end{equation}
and $\widehat a_0(L/2)$ is introduced   in \eqref{eq:def-hat-a}.
\end{lemma}
\begin{proof}
We express the function $w$ introduced in Lemma~\ref{lem:J1} as
\[w(x_2)=f(x_2)+\ii g(x_2),\quad f(x_2)=\cos\theta(x_2),\quad g(x_2)=\sin\theta(x_2),\]
 where $\theta(x_2)=2\gamma_0\arctan(2x_2/L)$. Then, 
\begin{multline*}w(x_2)\bigl(F(x_2)-\ii G(x_2)\bigr)\\
=\bigl(f(x_2)F(x_2)+g(x_2)G(x_2)\bigr)+\ii \bigl(g(x_2)F(x_2)-f(x_2)G(x_2)\bigr). \end{multline*}
Observing that the functions $f,F$ are even while the functions $g,G$ are odd, we infer from Lemma~\ref{lem:J1},
\begin{equation}\label{eq:J1*}
J_1=2\pi^{-1}h^2\ee^{-2\ii\pi\gamma_0}\int_{0}^{+\infty} \bigl(f(x_2)F(x_2)+g(x_2)G(x_2)\bigr)\dd x_2.
\end{equation}
We do the change of variable $r=\sqrt{x_2^2+(L/2)^2}$ and get
\[ 
\begin{gathered}
\int_{0}^{+\infty} g(x_2)G(x_2)\dd x_2=2\gamma_0\int_{L/2}^{+\infty} g\bigl(\sqrt{r^2-(L/2)^2}\,\bigr) r^{-1}|\psi(r)|^2\,\dd r,\\
\int_{0}^{+\infty}f(x_2)F(x_2)\dd x_2=L\int_{L/2}^{+\infty}\frac{f\bigl(\sqrt{r^2-(L/2)^2}\,\bigr)}{\sqrt{r^2-(L/2)^2}}\,\psi(r)\,\psi'(r)\,\dd r.
\end{gathered}
\]
Choose $R>L$ such that $d(R)> 2d(L/2)$, then by Proposition~\ref{prop:dec-psi},
\begin{multline*}
\int_{0}^{+\infty} g(x_2)G(x_2)\dd x_2\\=2\gamma_0\int_{L/2}^{R} g\bigl(\sqrt{r^2-(L/2)^2}\,\bigr) r^{-1}|\psi(r)|^2\,\dd r+o\bigl(\ee^{-2d(L/2)/h}\bigr), \end{multline*}
and
\begin{multline*}
\int_{0}^{+\infty}f(x_2)F(x_2)\dd x_2\\=L\int_{L/2}^{R}\frac{f\bigl(\sqrt{r^2-(L/2)^2}\,\bigr)}{\sqrt{r^2-(L/2)^2}}\,\psi(r)\,\psi'(r)\,\dd r+o\bigl(\ee^{-2d(L/2)/h}\bigr).
 \end{multline*}
We focus now on computing the integral involving $F$.
By the WKB approximation in Theorem~\ref{thm:T0}, we have, uniformly on $[L/2,R]$,  
\begin{equation*}
\psi(r)\, \psi'(r)= -h^{-2-e_0}\sqrt{|v(0)|}\, r^{2e_0}\bigl(|\widehat a_0(r)|^2+o(1)\bigr)\ee^{-2d(r)/h} \,. \end{equation*} 
Moreover, for $r\geq L/2$, we have
\[d(r)=d(L/2)+\Bigl(r-\frac{L}{2}\Bigr)\sqrt{|v(0)|}\,,\]
and for $r=L/2$, \[f\bigl(\sqrt{r^2-(L/2)^2}\,\bigr)=f(0)=1\,. \]

Consequently, with \[ M=\sqrt{2|v(0)|}(L/2)^{2e_0+\frac12}|\widehat a_0(L/2)|^2\,,\]
we have
\begin{multline*}
L\int_{L/2}^{+\infty}\frac{f\bigl(\sqrt{r^2-(L/2)^2}\,\bigr)}{\sqrt{r^2-(L/2)^2}}\,\psi(r)\,\psi'(r)\,\dd r\\
= -\bigl(1+o(1)\bigr)Mh^{-2-e_0} \ee^{-2d(L/2)/h}\int_{L/2}^{+\infty}\frac{1}{\sqrt{r-L/2}}\ee^{-2(r-\frac{L}{2})\sqrt{|v(0)|}/h}\,\dd r\\
=  -\frac{\pi^{\frac12}(L/2)^{2e_0+\frac12}|\widehat a_0(L/2)|^2}{|v(0)|^{\frac14}} h^{-\frac32-e_0}\, \ee^{-2d(L/2)/h} + o\bigl(h^{-\frac32-e_0}\, \ee^{-2d(L/2)/h}\bigr).
\end{multline*}
We handle the integral involving $G$ in a similar manner. Firstly, by  Theorem~\ref{thm:T0}, we have
\[
|\psi(r)|^2= -h^{-1-e_0}\, r^{2e_0}\bigl(|\widehat a_0(r)|^2+o(1)\bigr)\ee^{-2d(r)/h}\mbox{ uniformly on }[L/2,R]).  \]
Secondly, we use this approximation to write 
\[\int_{L/2}^{R} g\bigl(\sqrt{r^2-(L/2)^2}\,\bigr)\, r^{-1}|\psi(r)|^2\,\dd r=\cO\bigl(h^{-e_0}\ee^{-2d(L/2)/h}\bigr), \]
and consequently
\[\int_{0}^{+\infty} g(x_2)G(x_2)\dd x_2 =o\bigl(h^{-\frac32-e_0}\ee^{-2d(L/2)/h}\bigr).\]
\end{proof}

\begin{proof}[Proof of Theorem~\ref{thm:main}]
Collect \eqref{eq:J1*} and \eqref{eq:splitting}, and notice that $2S_\varepsilon>2d(L/2)$ by \eqref{eq:def-S-ep}. Then, we have the asymptotics  in Theorem~\ref{thm:main}. The uniform convergence for $e_0\in[0,\delta_0)$ holds as explained in Remark~\ref{rem:T0-unif} and i) in Subsection~\ref{subsec:dw-op}.
\end{proof}
\begin{proof}[The spectral gap for the third and fourth excited states]~~

To estimate $\lambda_4(h,\alpha)-\lambda_3(h,\alpha)$, we repeat  the analysis leading to the proof of Theorem~\ref{thm:main}. Let us explain the due adjustments: We use the abstract result from \cite{HKS} but for the operator $\cT=\sH-\lambda_2^\sw$ and the Hilbert space $L^2(\R^2)\ominus F$, where $\lambda_2^\sw$ is the second eigenvalue of the single well operator $\sH_0$, and 
$F=\oplus_{i=1}^2\mathrm{Ker}\bigl(\sH-\lambda_i(h,\alpha)\bigr)$ is the span of the ground and second excited states of the double well operator. After separation of variables, as we did in Subsection~\ref{subsec:sep-var},  we encounter 
\[\hat e_0=\min_{\substack{m\in\Z\\ m\not=m_*}}|m-\alpha/h|, \]
where $m_*$ is introduced in \eqref{eq:def-m} (excluding $m_*$ amounts to excluding the angular momentum that corresponds to the ground state of the single well operator). We then build the quasi-modes using the positive normalized ground state of the $1$D operator 
\[\sT_{\hat e_0}=-\frac{\dd^2}{\dd r^2}-\frac12\frac{\dd}{\dd r}+v+\frac{\hat e_0^2}{r^2}.\]
For $0<e_0<\frac12$, we get $\hat e_0=1-e_0\in(0,\frac12)$, and we notice that $\lambda_2(\sT_{e_0})>\lambda_1(\sT_{\hat e_0})$, thanks to \eqref{eq:2.2}. Then,  \eqref{eq:lambda4-3} results from the same computations leading to Theorem~\ref{thm:main}, but with $\hat e_0=1-e_0$ instead of $e_0$. 

Note that this argument requires that $e_0\not=0$. However, if $e_0=0$ the accurate estimate of $\lambda_4(h,\alpha)-\lambda_3(h,\alpha)$  was carried out in \cite{HSj}, because $\sH$ is unitarily equivalent to  the flux free operator.
\end{proof}

\begin{remark}[Symmetry relations of eigenfunctions]\label{rem5.3}

Since  $\sL_\Psi$ commutes with $\sH$ and $L_\Psi^2=\sI$, we  can always find an orthonormal basis of eigenfunctions $u_j$ of $\sH$ and $k_j \in\{-1,1\}$ such that 
\[\sL_\Psi u_j=k_j u_j.\]
\end{remark}

\begin{remark}[Locally radial wells]~~~
We can  handle the situation when the potential $V$ is a smooth function that satisfies
\[\begin{cases}
V(x_1,x_2)=V(-x_1,x_2)\quad (\xv=(x_1,x_2)\in\R^2),\\
V_0:=\min V=\{V(\xv_\ell),V(\xv_r)\},\\
V\mbox{ is radial in a neighborhood of $\xv_\ell$ (respectively $\xv_r$),}\\
V''(\xv_\ell)>0\mbox{ and }V''(\xv_r)>0.
\end{cases} \]
Denote by $d(\xv,\mathsf{y})$ the Agmon distance between $\xv$ and $\yv$ with respect to the metric $\max(V(\xv)-V_0,0)\dd \xv^2$, and let $S(v,L)= d(\xv_\ell,\xv_r) $. Then, assuming that
there is a unique minimal geodesic connecting $\xv_\ell$ and $\xv_r$ (which is non degenerate in the sense of \cite{HSj}), we can prove a similar statement to Theorem~\ref{thm:main}. In fact, we have a radially symmetric WKB approximation  of the  single well ground state near the minimum of the well, that we extend along the geodesic connecting $\xv_\ell$ and $\xv_r$. We skip the technical details and refer to \cite{HSj} (see also \cite[Thm.~4.1]{MR}). 
\end{remark}

\section{Half-integer flux}\label{sec:interaction*}

We consider $e_0=1/2$ and suppose that $h$ varies in $\cJ_\alpha(e_0)$, so that $\alpha/h\in \Z+\frac12$. We recall by \eqref{eq:def-gs-AB*} that the ground state energy $\lambda$ of the  single well operator has multiplicity two and the corresponding  linearly independent and normalized ground states
\[\phi_1^\sw=\phi^\sw,\quad \phi_2^\sw=\hat\phi^\sw:=\ee^{\ii\theta}\phi^\sw, \]
where  
\begin{equation}\label{eq:phi-sw*}
\phi^\sw=\pi^{-\frac12}\psi(r)\ee^{\ii m_*\theta},\quad m_*=m_*(\alpha/h):=\frac{\alpha}{h}-\frac12,
\end{equation}
and $\psi$ is the positive normalized ground state of $\sT$. For instance, observe that in this case (see \eqref{eq:def-gam0}),
\[\gamma_0:=m_*-\alpha/h=-1/2.\]
We recall that $\eta(z)=\theta$ if $(r,\theta)\in\R_+\times[0,2\pi)$ are the polar coordinates of $z\not=0$,  and for simplicity of notation, we write  $\sH$ for the operator $\sH_{\alpha}$ introduced in \eqref{eq:def-op}. We will encounter the operators
\[\begin{aligned}
\sH_0^\ell&=\bigl(-\ii h\nabla-\Fb_\alpha(\xv-\xv_\ell)\bigr)^2+\vb(\xv-\xv_\ell),\\
 \sH_0^r&=\bigl(-\ii h\nabla-\Fb_\alpha(\xv-\xv_r)\bigr)^2+\vb(\xv-\xv_r),
 \end{aligned}\] 
 where $\Fb_\alpha$ was introduced in \eqref{eq:def-F}.

\subsection{The quasi-modes}~\medskip

We recall that $\varepsilon>0$ is small but fixed, and  that there is a corresponding constant  $S_\varepsilon$ satisfying \eqref{eq:def-S-ep}. 

We recall the two functions defined on $\cS:=\{\xv=(x_1,x_2)\colon |x_1|<L/2\}$ as
\[\phi_\ell(\xv)=\phi^\sw(\xv -\xv_\ell)\ee^{\ii \alpha\eta(\xv-\xv_r)/h},\quad \phi_r(\xv)=( \bar\sL_\Phi \phi_\ell)(\xv), \]
where $ \bar\sL_\Phi$ is the operator introduced in \eqref{eq:reflection}. It will be convenient to introduce the two functions
\[f_\ell(x)=\chi_\ell(\xv)\phi^\sw(\xv -\xv_\ell),\quad f_r(\xv)=(\bar\sL_\Phi f_\ell)(\xv)= -\chi_r(\xv)\phi^\sw(\xv -\xv_r).\]
 Let the quasi-modes $u_\ell$ and $u_r$ be as in \eqref{eq:qm-ell} and \eqref{eq:qm-r} respectively. They  can be expressed as
\begin{equation}\label{eq:ur-ul}
u_\ell(\xv )=f_\ell(\xv)\ee^{\ii \alpha\eta(\xv-\xv_r)/h},\quad u_r(\xv )=( \bar\sL_\Phi u_\ell)(\xv),
\end{equation}
and we recall the  useful identity,
\begin{equation}\label{eq:red-1well}
\sH u_\ell=\ee^{\ii \alpha\eta(\xv-\xv_r)/h}\sH_0^\ell f_\ell. 
\end{equation}
We also introduce the quasi-modes (with support in $\cS$),
\begin{equation}\label{eq:hat-ur-ul}
\hat u_\ell=\ee^{\ii \eta(\xv -\xv_\ell)} u_\ell,\quad \hat u_r= \bar\sL_\Phi\hat u_\ell=-\ee^{\ii\eta(\xv -\xv_r)}u_r.
\end{equation}
It is convenient to express $\hat u_\ell$ via the anti-linear operators $\sK_\ell$ defined as follows (the definition makes sense because $2\alpha/h\in\Z$)
\[u\mapsto \sK_\ell=\ee^{2\ii\alpha\eta(\xv-\xv_\ell)/h}\,\overline u.  \]
In fact, by \eqref{eq:phi-sw*} and \eqref{eq:ur-ul},
\begin{equation}\label{eq:exp-hat-u-ell}
\hat u_\ell=\ee^{\ii \alpha\eta(\xv-\xv_r)/h}\sK_\ell f_\ell=\ee^{2\ii\alpha\eta(\xv-\xv_r)/h}\ee^{2\ii\alpha\eta(\xv-\xv_\ell)/h}\,\overline{u_\ell},
\end{equation}
from which we derive the following formula
\[\hat u_r(\xv)=\overline{\hat u_\ell(\Phi(\xv))}=\ee^{-2\ii\alpha\eta(\Phi(\xv)-\xv_r)/h}\ee^{-2\ii\alpha\eta(\Phi(\xv)-\xv_\ell)/h}\overline{u_r(\xv)}. \]
In light of \eqref{eq:def-op-K*}, \eqref{eq:red-1well} and \eqref{eq:hat-ur-ul}, we have
\[\begin{aligned}
\sH\hat u_\ell&=\ee^{\ii \alpha\eta(\xv-\xv_r)/h}\sK_\ell\sH_0^\ell f_\ell\\
&=\ee^{2\ii\alpha\eta(\xv-\xv_r)/h}\ee^{2\ii\alpha\eta(\xv-\xv_\ell)/h}\,\overline{\sH u_\ell}\,.
\end{aligned} \]
Thanks to \eqref{eq:symm-eta} and the fact that $2\alpha/h\in\Z$, we deduce the following identities
\[
\langle \hat u_\ell,\hat u_r\rangle=\overline{\langle u_\ell,u_r\rangle},\quad \langle \sH\hat u_\ell,\hat u_r\rangle=\overline{\langle \sH u_\ell,u_r\rangle},\quad \langle\sH\hat u_\ell,\hat u_\ell\rangle=\overline{\langle \sH u_\ell,u_\ell\rangle}.\]
Moreover,  by \eqref{eq:symm-H} and \eqref{eq:ur-ul}, we have that $\langle\sH u_\ell,u_\ell\rangle=\langle \sH u_r,u_r\rangle$. Consequently,
\[
\begin{gathered}
\langle \cT u_\ell,u_\ell\rangle=\langle \cT u_r,u_r\rangle=J_0\\
\langle \cT u_\ell,u_r\rangle=\overline{\langle \cT u_r,u_\ell\rangle}=J_1\\
\langle \cT \hat u_\ell,\hat u_\ell\rangle=\langle \cT \hat u_r,\hat u_r\rangle=J_0\\
\langle \cT \hat u_\ell,\hat u_r\rangle=\overline{\langle \cT \hat u_r,\hat u_\ell\rangle}=
\overline{J_1}
\end{gathered}
\]
where $J_0$ and $J_1$ are as in \eqref{eq:J1,2}, and $\cT=\sH-\lambda$. In particular, we have an estimate of $J_0$ and an accurate approximation of $J_1$.

As for the other interaction terms, we also have by \eqref{eq:ur-ul} and \eqref{eq:hat-ur-ul},
\[
\begin{gathered}
\langle \cT \hat u_\ell,u_\ell\rangle=\overline{\langle \cT \hat u_r,u_r\rangle}=\hat J_0\\
\langle \cT \hat u_\ell,u_r\rangle=\overline{\langle \cT u_r,\hat u_\ell\rangle}=\hat J_1\\
\langle \cT u_\ell,\hat u_\ell\rangle=\overline{\langle \cT u_r,\hat u_r\rangle}=\overline{\hat J_0}\\
\langle \cT u_\ell,\hat u_r\rangle=\overline{\langle \cT \hat u_r,u_\ell\rangle}=\hat J_1
\end{gathered}
\]
and by  the same reasoning   used to estimate  $J_0$ and $J_1$, we have
\[ |\hat J_0|=\cO(\ee^{-S_\varepsilon/h}),\quad |\hat J_1|=\cO(\ee^{-S_\varepsilon/h}).\]
Consequently, there is a basis of the vector space spanned by the eigenfunctions corresponding to the eigenvalues $\{\lambda_j(\sH)\colon 1\leq j\leq 4\}$ such that, the matrix $\sW$  of $\cT=\sH-\lambda$ in this basis can be decomposed as follows (see \cite[Prop.~2.4]{HKS})
\[ \sW=\sU+\sR\]
where
\[\sR=\mathcal O(\ee^{-2S_\varepsilon/h}),\]
 and the (interaction) matrix $\sU$ is
\[\sU=
\begin{pmatrix}
\langle \cT u_\ell,u_\ell\rangle&\langle \cT u_\ell, u_r\rangle&\langle \cT u_\ell, \hat u_\ell\rangle&\langle \cT u_\ell, \hat u_r\rangle\\
\langle \cT u_r, u_\ell\rangle&\langle \cT u_r, u_r\rangle&\langle \cT u_r, \hat u_\ell\rangle&\langle \cT u_r, \hat u_r\rangle\\
\langle \cT \hat u_\ell, u_\ell\rangle&\langle \cT \hat u_\ell, u_r\rangle&\langle \cT \hat u_\ell, \hat u_\ell\rangle&\langle \cT \hat u_\ell, \hat u_r\rangle\\
\langle \cT \hat u_r, u_\ell\rangle&\langle \cT \hat u_r, u_r\rangle&\langle \cT \hat u_r, \hat u_\ell\rangle&\langle \cT \hat u_r,\hat u_r\rangle
\end{pmatrix}.
\]
In terms of the interaction coefficients $J_0,J_1,\hat J_0,\hat  J_1,$ introduced earlier, the matrix $\sU$ reads as
\[ \sU=
\begin{pmatrix}
J_0&J_1&\overline{\hat J_0}&\hat J_1\\
\overline{J_1}&J_0&\overline{\hat J_1}&\hat J_0\\
\hat J_0&\hat J_1&J_0&\overline{J_1}\\
\overline{\hat J_1}&\overline{\hat J_0}&J_1&J_0
\end{pmatrix}.\]

\subsection{The interaction terms}

We now discuss the calculation of the interaction terms $J_1,\hat J_1,J_0$ and $\hat J_0$. 

\subsubsection*{Calculation of $J_1$}

Starting with $J_1$,  we know by \eqref{eq:phi-sw*} that $\gamma_0=-\frac12$ in \eqref{eq:def-gam0}. It then follows from Lemma~\ref{lem:J1*},
\begin{equation}\label{eq:value-J1}
J_1= |J_1|=\frac12C(L,v,\mbox{$\frac12$}) \ee^{-2d(L/2)/h}+o\Bigr(\ee^{-2d(L/2)/h} \Bigr).
\end{equation}

\subsubsection*{Calculation of $\hat J_1$}
We introduce the two functions
\begin{equation}\label{eq:def-phi-hat}
\hat\phi_\ell(\xv )=\ee^{\ii \eta(\xv -\xv_\ell)}\phi_\ell(\xv ),\quad \hat\phi_r(\xv )=\ee^{\ii \eta(\xv -\xv_r)}\phi_r(\xv ),
 \end{equation}
and notice that by \eqref{eq:ur-ul} and \eqref{eq:hat-ur-ul},
\[ \hat u_\ell=\chi_\ell\hat\phi_\ell,\quad \hat u_r=\chi_r\hat\phi_r,\]
where $\chi_\ell$ and $\chi_r$ are the cut-off functions appearing in \eqref{eq:qm-ell} and \eqref{eq:qm-r} respectively. The conditions on the support of the cut-off functions $\chi_\ell$ and $\chi_r$ yield 
\[\hat J_1=\langle[\sH,\chi_\ell]\phi_\ell,\chi_r\hat\phi_r\rangle.\]
Then we do an integration by parts with respect to the variable $x_1$  and get
\[\begin{aligned}
\hat J_1&=\ii h \int_{\R} \left(P_1\phi_\ell\,\overline{\hat\phi_r}+\phi_\ell\,\overline{P_1\hat\phi_r}\right)(0,x_2)\,\dd  x_2\\
&=h^2\int_{\R}\bigl(\partial_1\phi_\ell\,\overline{\hat\phi_r}-\phi_\ell\,\overline{\partial_1\hat\phi_r}\bigr)(0,x_2)\, \dd  x_2-2\ii h\int_{\R}(A_1\phi_\ell\overline{\hat\phi_r})(0,x_2)\,\dd  x_2.
\end{aligned}\]
Knowing that $\hat\phi_r=\ee^{\ii\eta(\xv -\xv_r)}\phi_r$, we have
\begin{multline*}
\hat J_1=h^2\int_{\R}\ee^{-\ii\eta(\xv -\xv_r)}\bigl(\partial_1\phi_\ell\,\overline{\phi_r}-\phi_\ell\,\overline{\partial_1\phi_r}\bigr)(0,x_2)\, \dd  x_2\\
-\ii h\int_\R\ee^{-\ii\eta(\xv -\xv_r)}\Bigl(2A_1- h \partial_1\eta(\xv -\xv_r)\Bigr)\phi_\ell\overline{\phi_r}(0,x_2)\,\dd  x_2.
\end{multline*}
The computations are then identical to those carried out in Lemma~\ref{lem:J1}. Noticing that $\gamma_0=-\frac12$, and using \eqref{eq:relation-eta} (and the subsequent computations of the integrand), we have
\[\ee^{-\ii\eta(\xv -\xv_r)}\big|_{x_1=0}=-\ee^{\ii\arctan(2x_2/L)}= -1/w(x_2),\]
hence
\[
\hat J_1= -\pi^{-1}h^2\int_\R\bigl(F(x_2)-\ii G(x_2)\bigr)\dd  x_2\\
-h^2 \int_\R\frac{x_2\ee^{-\ii\eta(\xv -\xv_r)}}{(x_2^2+(L/2)^2)}\phi_\ell\overline{\phi_r}(0,x_2)\dd  x_2,
\]
where $F$ and $G$ are introduced in \eqref{eq:def-F,G}, with $\gamma_0=-1/2$. 

Using \eqref{eq:qm-p'}, we get
\[
\hat J_1= -\pi^{-1}h^2 \int_\R\bigl(F(x_2)-\ii G(x_2)\bigr)\dd  x_2
-\pi^{-1}h^2 \int_\R G(x_2)\ee^{\ii\arctan(2x_2/L)}\dd  x_2.
\]
Recall that $F$ is even and $G$ is odd. Omitting the odd functions appearing in the integrand, the expression of $\hat J_1$ becomes
\[
\hat J_1=-2\pi^{-1}h^2 \int_0^{+\infty} F(x_2)\dd  x_2-2\pi^{-1}h^2 \int_0^{+\infty} G(x_2)\ee^{\ii\arctan(2x_2/L)}\dd  x_2.
\]
We proceed by using the WKB approximation and the Laplace method, as in Lemma~\ref{lem:J1*}. Eventually we get
\begin{equation}\label{eq:hat-J1}
\hat J_1=-\frac12C(L,v,\mbox{$\frac12$})\ee^{-2d(L/2)/h}+o\bigl(\ee^{-2d(L/2)/h}\bigr)=-J_1+o(J_1).
\end{equation}
\subsubsection*{Estimating $J_0$ and $\hat J_0$}
We will prove that 
\begin{equation}\label{eq:comp-J0}
J_0=0\quad\mbox{and}\quad \hat J_0=o(J_1).
\end{equation}
Similarly as we did for $J_1$ in \eqref{eq:J1} and  \eqref{eq:J1*}, we have
\[\begin{aligned}
J_0&=\langle [\sH,\chi_\ell]\phi_\ell,\chi_\ell\phi_\ell\rangle\\
&=h^2\int_{\R}\bigl(\partial_1\phi_\ell\,\overline{\phi_\ell}-\phi_\ell\,\overline{\partial_1\phi_\ell}\bigr)(0,x_2)\, \dd  x_2-2\ii h\int_{\R}(A_1\phi_\ell\overline{\phi_\ell})(0,x_2)\,\dd  x_2\\
&=-\pi^{-1}h^2\ii\int_\R G(x_2)\dd  x_2=0,
\end{aligned} \]
where $G$ is the odd function introduced in \eqref{eq:def-F,G}, with $\gamma_0=-1/2$.

We do a similar calculation for $\hat J_0$ and write
\[\begin{aligned}
\hat J_0&=\langle [\sH,\chi_\ell]\hat\phi_\ell,\chi_\ell\phi_\ell\rangle\\
&=h^2\int_{\R}\bigl(\partial_1\hat\phi_\ell\,\overline{\phi_\ell}-\hat\phi_\ell\,\overline{\partial_1\phi_\ell}\bigr)(0,x_2)\, \dd  x_2-2\ii h\int_{\R}(A_1\hat \phi_\ell\overline{\phi_\ell})(0,x_2)\,\dd  x_2.
\end{aligned} \]
Using that $\hat\phi_\ell=\ee^{\ii\eta(\xv -\xv_\ell)}\phi_\ell$, then up to the prefactor $\ee^{\ii\eta(\xv -\xv_\ell)}$, the integrand becomes 
\[\Bigl(h^2\bigl(\partial_1\phi_\ell\,\overline{\phi_\ell}-\phi_\ell\,\overline{\partial_1\phi_\ell}\bigr)-2\ii h(A_1|\phi_\ell|^2-h^2|\phi_\ell|^2\partial_1\eta(\xv -\xv_\ell)\Bigr)(0,x_2). \]
Consequently, as we did in the computations leading to \eqref{eq:J1*}, we get
\begin{multline*}\hat J_0=-\pi^{-1}h^2\ii\int_\R \ee^{\ii\eta(\xv -\xv_\ell)}G(x_2)\dd  x_2\\
-\pi^{-1}h^2\ii\int_{\R}
\partial_1\eta(\xv -\xv_\ell)\big|_{x_1=0}\Bigl|\psi\Bigl(\sqrt{x_2^2+(L/2)^2}\Bigr)\Bigr|^2\dd  x_2. 
\end{multline*}
Using Theorem~\ref{thm:T0} and Proposition~\ref{prop:dec-psi}, we get
\[\hat J_0=o\bigl(\ee^{-2d(L/2)/h}\bigr)\]
and thanks to \eqref{eq:comp-J1}, this yields that $\hat J_0=o(J_1)$.

\subsection{Eigenvalues of the interaction matrix}

We now have explicit formulas for all the entries of the interaction matrix, thanks to \eqref{eq:value-J1}, \eqref{eq:hat-J1}  and  \eqref{eq:comp-J0}. So we can express the interaction matrix as
\[\sU=|J_1|\bigl(\hat\sU+\hat\sR\bigr), \]
where
\[\hat\sU=\begin{pmatrix}
0&1&0&-1\\
1&0&-1&0\\
0&-1&0&1\\
-1&0&1&0
\end{pmatrix}\]
and $\hat\sR=o(1)$ as $h\to0$. 

The eigenvalues of $\hat U$ are
\[
\lambda_1(\hat U)=-2,\quad  \lambda_2(\hat U)=\lambda_3(\hat U)=0,\quad
\lambda_4(\hat U)=2\,, 
\]
corresponding with the following basis of eigenvectors
\[
 \Bigl\{\bp_1=\frac 12\begin{pmatrix}-1\\1\\1\\-1\end{pmatrix},~   \bp_2=\frac{1}{\sqrt{2}}
\begin{pmatrix}
 1\\0\\1\\0\end{pmatrix},~ \bp_3=\frac{1}{\sqrt{2}} \begin{pmatrix}0\\1\\0\\1\end{pmatrix},~\bp_4= \frac 12\begin{pmatrix}1\\1\\-1\\-1\end{pmatrix}\Bigr\} \,.
\]

Returning back to the eigenvalues of $\sH$, which are the eigenvalues of the matrix $\sW$ of $\cT=\sH-\lambda$, we get
\begin{equation}\label{eq:ev-half-integer}
\lambda_k(\sH)=\lambda+|J_1|\lambda_k(\hat\sU) +o(|J_1|),\quad (1\leq k\leq 4).
 \end{equation}
 This yields the conclusion in Theorem~\ref{thm:main*}.

\begin{remark}\label{rem:ef-mult}
There is an orthonormal  basis of  eigenvectors  $\{\bq_1,\bq_2,\bq_3,\bq_4\}$ for the matrix $\sW$ such that
\[ \begin{gathered}
|\bq_1-\bp_1|=o(1),\quad |\bq_4-\bp_4|=o(1),\\
\mathrm{dist}\bigl(\bq_i,\mathrm{Span}(\bp_2,\bp_3)\bigr)=o(1)\mbox{ for }i=2,3.
\end{gathered}\]
With
\[\left\{ \begin{gathered}
\hat\bp_1=-\frac{1}{2}(u_\ell-u_r)+\frac{1}{2}(\hat u_\ell-\hat u_r),\\
\hat\bp_2=\frac{1}{\sqrt{2}}(u_\ell+\hat u_\ell),\\
\hat\bp_3=\frac{1}{\sqrt{2}}(u_r+\hat u_r),\\
\hat\bp_4=\frac{1}{2}(u_\ell+u_r)-\frac{1}{2}(\hat u_\ell+\hat u_r),
\end{gathered}\right.\]
there is a corresponding basis of orthonormal eigenfunctions  of $\sH$
\[ \{\hat\bq_1,\hat\bq_2,\hat\bq_3,\hat\bq_4\}\]
such that
\[
\begin{gathered}
\|\hat\bq_1-\hat\bp_1\|_{L^2(\R^2)}=o(1),\quad \|\hat\bq_4-\hat\bp_4\|_{L^2(\R^2)}=o(1),\\
\mathrm{dist}_{L^2(\R^2)}\bigl(\hat\bq_i,\mathrm{Span}(\hat\bp_2,\hat\bp_3)\bigr)=o(1).
\end{gathered}\]
\end{remark}

\subsection{Analysis of the multiplicity}\label{subsec:mult}

\subsubsection*{Using an anti-linear operator}~\\
Let us start by noticing that
\[\begin{gathered}
\Fb(\xv-\xv_\ell)=\nabla\eta(\xv-\xv_\ell)\mbox{ on }\R^2\setminus\{\xv_\ell\},\\
\Fb(\xv-\xv_r)=\nabla\eta(\xv-\xv_r)\mbox{ on }\R^2\setminus\{\xv_r\}.
\end{gathered}\]
We introduce the function defined on $\R^2\setminus\{\xv_\ell,\xv_r\}$ as
\[\zeta(\xv)=2(\alpha/h)\bigl(\eta(\xv-\xv_\ell)+\eta(\xv-\xv_r) \bigr),\]
and  the vector potential in \eqref{eq:def-A} satisfies
\[2\Ab=h\nabla\zeta \mbox{ on }\R^2\setminus\{\xv_\ell,\xv_r\}.\]
Since $2\alpha/h$ is an integer, we can introduce the anti-linear operator
\[ \sK^{\dw}:L^2(\R^2)\ni u\mapsto \ee^{\ii \zeta} \bar u\in L^2(\R^2),\]
and we notice that
\[   (\sK^{\dw})^2=\sI,\quad \sH\sK^{\dw}=\sK^{\dw}\sH.\]
To calculate $\sK^{\dw}\hat\bq_2$, we first calculate 
\[
\sK^{\dw}\hat\bp_2=\frac1{\sqrt{2}}\ee^{\ii \zeta}\bigl(\overline{u_\ell}+\overline{\hat u_\ell}\bigr)=\hat\bp_2,
\]
where we used \eqref{eq:exp-hat-u-ell} to verify that $\sK^{\dw}u_\ell=\hat u_\ell$, and consequently $\sK^{\dw}\hat u_\ell=u_\ell$. We can also verify that $\sK^{\dw}u_r=\hat u_r$ and $\sK^{\dw}\hat u_r= u_r$.  This leads to
\[
\sK^{\dw}\hat\bp_3=\frac1{\sqrt{2}}\ee^{\ii \zeta}\bigl(\overline{u_\ell}+\overline{\hat u_\ell}\bigr)=\hat\bp_3,
\]
Similar calculations lead to
\[
\sK^{\dw}\hat\bp_1=- \hat\bp_1 , \sK^{\dw}\hat\bp_4=- \hat\bp_4.
\]
In conclusion, $\sK^{\dw}$ leaves invariant each eigenvector. This cannot explain the occurrence of  multiplicity.

\subsubsection*{Using the symmetry operator $L_\Psi$.}~\\
We recall that
\[\sL_\Psi ^2=\sI,\quad \sH\sL_\Psi =\sL_\Psi \sH.\]
If $\xv=(x_1,x_2)$ with $x_1>-L/2$, then both $\eta(\Phi(\xv-\xv_\ell))$ and $\eta(-\xv-\xv_r)$ are in $(0,2\pi)$, and we have  $\eta(-\xv-\xv_r)=2\pi-\eta(\Phi(\xv-\xv_\ell))$; moreover, $\Phi(\xv-\xv_\ell)=\Phi(\xv)-\xv_r$. A similar relation exists between $\eta(-\xv-\xv_\ell)$ and $\eta(\Phi(\xv-\xv_r))$.  Consequently, by \eqref{eq:qm-ell} and \eqref{eq:qm-r},
\[(\bar\sL_\Phi u_\ell)(\xv)=u_\ell(-\xv)=\overline{u_\ell(\Phi(\xv))}=u_r(\xv). \]
Observing  that  $\eta(-\mathsf{y})=\eta(\mathsf{y})+\pi$ modulo $2\pi$, it results from \eqref{eq:hat-ur-ul}, 
\[u_\ell(\xv)+\hat u_\ell(\xv)=\bigl(1+\ee^{\ii\eta(\xv-\xv_\ell)} \bigr)u_\ell(\xv),\quad u_r(\xv)+\hat u_r(\xv)=\bigl(1-\ee^{\ii\eta(\xv-\xv_r)} \bigr)u_r(\xv),\]
and consequently
\[\begin{aligned}(\sL_\Psi u_\ell)(\xv)+( \sL_\Psi\hat u_\ell)(\xv)
&=\bigl(1-\ee^{\ii\eta(\xv-\xv_r)}\bigr)u_\ell(-\xv)\\
&=u_r(\xv)+\hat u_r(\xv).
\end{aligned} \]
Thus, we verified that 
\[\sL_\Psi\hat\bp_2=\hat \bp_3,\quad  \sL_\Psi\hat\bp_3=\hat \bp_2. \]
 In a similar fashion, we can verify that
\[\sL_\Psi \hat \bp_1=-\hat\bp_1,\quad \sL_\Psi \hat \bp_4=\hat\bp_4\,. \]
 Therefore, using Remark \ref{rem5.3} and the relations satisfied by the $\hat \bp_j$, there is a corresponding orthonormal basis $\{\tilde\bq_1,\tilde\bq_2,\tilde\bq_3,\tilde\bq_4 \}$ of  joint eigenfunctions  of $\sH$  (restricted to the eigenspace corresponding to the four first eigenvalues) and $\sL_\Psi$
 such that
\[
\begin{gathered}
 \sH \tilde\bq_1= \lambda_1(\sH)  \tilde\bq_1\,,\,  \sH  \tilde\bq_2 =\lambda_2(\sH)  \tilde\bq_2\,,\,  \sH  \tilde\bq_3=\lambda_3(\sH)  \tilde \bq_3\, ,\,  \sH \tilde\bq_4 = \lambda_4 (\sH)  \tilde \bq_4 \\
\sL_\Psi \tilde\bq_1=-  \tilde\bq_1\,,\,  \sL_\Psi \tilde\bq_2 =\pm \tilde\bq_2\,,\,  \sL_\Psi \tilde\bq_3=\mp \tilde \bq_3\, ,\,  \sL_\Psi \tilde\bq_4 = \tilde \bq_4 \\
\sK^{\dw} \tilde \bq_1=-\tilde\bq_1\,,\,\sK^{\dw} \tilde\bq_2=\tilde\bq_2\,,\, \sK^{\dw} \tilde\bq_3=\tilde\bq_3\,,\, \sK^{\dw} \tilde\bq_4=-\tilde\bq_4\,\\
\|\tilde\bq_1-\hat\bp_1\|_{L^2(\R^2)}=o(1),\quad \|\tilde\bq_4-\hat\bp_4\|_{L^2(\R^2)}=o(1),\\
\|\tilde\bq_2 -(\hat\bp_2\pm \hat\bp_3)\|_{L^2(\R^2)}=o(1),\quad \|\tilde\bq_3-(\hat\bp_2 \mp \hat\bp_3)\|_{L^2(\R^2)}=o(1) \,.
\end{gathered}
\]

\subsection*{Acknowledgments}

\footnotesize{AK is partially supported by CUHK-SZ grant no. UDF01003322 and UF02003322}.

\end{document}